\documentclass[12pt,a4paper]{amsart}
\usepackage{newunicodechar}
\newunicodechar{​}{} %
\usepackage[top=35mm, bottom=35mm, left=30mm, right=30mm]{geometry}
\usepackage{mathptmx}
\usepackage{mathrsfs}
\usepackage{amscd}
\usepackage{verbatim}
\usepackage{url}
\usepackage[all]{xy}
\usepackage{color}
\usepackage[colorlinks=true,citecolor=blue, backref]{hyperref}%
\usepackage{amsmath}
\usepackage{amssymb}
\usepackage{array}

\theoremstyle{plain}

\newtheorem{thma}{Theorem}

\newtheorem{thm}{Theorem}[section]
\newtheorem{lem}[thm]{Lemma}
\newtheorem{prop}[thm]{Proposition}
\newtheorem{cor}[thm]{Corollary}

\theoremstyle{definition}
\newtheorem{defn}[thm]{Definition}
\newtheorem{rem}[thm]{Remark}
\newtheorem{ques}{Question}

\newtheorem{ex}{Example}

\newcommand{\Z}{{\mathbb{Z}}}
\newcommand{\N}{\mathbb{N}}
\newcommand{\F}{\mathcal F}

\newcommand{\ep}{\varepsilon}
\newcommand{\ra}{\rightarrow}
\newcommand{\lra}{\rightarrow}
\newcommand{\T}{\mathbb{T}}

\DeclareMathOperator{\diam}{diam}

\def \id {{\rm id}}
\def\O {\mathcal O}

\begin{document}
\title[On subsets of integers having dense orbits]{On subsets of integers having dense orbits}

\author{Zhuowen Guo}
\address{$^{1,2,4}$School of Mathematical Sciences, University of Science and Technology of China, Hefei, Anhui, 230026, P.R. China}
\email{guozw0920@mail.ustc.edu.cn}
\author{Jiahao Qiu}
\email{qiujh@mail.ustc.edu.cn}
\author{Hui Xu}
\address{Department of Mathematics, Shanghai Normal University, Shanghai, 200234, CHINA}
\email{huixu@shnu.edu.cn}
\author{Xiangdong Ye}
\email{yexd@ustc.edu.cn}


\thanks{This research is supported by National Key R$\&$D Program of China (No. 2024YFA1013601, 2024YFA1013600), USTC Research Funds of the Double First-Class Initiative (YD0010002009), and National Natural Science Foundation of China (12031019, 12371196, 12401243, 12426201, 12471188).}

\maketitle

\begin{abstract}
Let $A\subset \mathbb{N}$. We say $A$ is an $R$-sequence for a given minimal system $(Y,S)$ if there is $y\in Y$ such that $\{S^ny:n\in A\}$ is dense in $Y$. Richter asked if $A$ is an $R$-sequence for all minimal equicontinuous systems implies that $A$ is an $R$-sequence for all minimal systems. In this paper, we investigate this question and related issues within the framework of totally minimal systems, including a characterization of transitive systems that are disjoint from all totally minimal systems.

A dynamical system is scattering (resp. weakly scattering) if its product with any minimal (resp. minimal and equicontinuous) system is transitive. It turns out that $(X,T)$ is scattering if and only if for any transitive point $x\in X$ and any minimal system $(Y,S)$ there is $y\in Y$ such that the orbit  of $(x,y)$ is dense in $X\times Y$ if and only if for each transitive point $x\in X$ and any non-empty open subset $U$ of $X$, $\{n\in \mathbb{N}:T^nx\in U\}$ is an $R$-sequence. By combining this result with earlier work of Huang and Ye, we deduce that if scattering and weak scattering are distinct properties, then both Richter's question and Katznelson's question admit negative answers.
\end{abstract}
\section{Introduction}

By a {\it topological dynamical system} (t.d.s. for short) we mean a compact metric space $X$ together with a continuous onto map $T:X\ra X$.

\subsection{Motivations}
Recurrence is one of the most fundamental concepts in dynamics, with wide-ranging applications across mathematics including combinatorial number theory.

\medskip
For a t.d.s. $(X,T)$, $x\in X$ is a {\it recurrent point} if there is a sequence $\{n_i\}_{i=1}^\infty$ in $ \N$ such that
$T^{n_i}x\ra x$. The well-known theorem of Birkhoff asserts that there is at least one recurrent point for each t.d.s.

A central open problem in the field connecting dynamics, combinatorial number theory, and group theory is the {\it Katznelson's question}. To frame it, we need to introduce some definitions.


We say that $A\subset \N$ is a {\it Birkhoff (resp. Bohr) recurrence sequence}, denoted by $A\in \F_{Bir}$ (resp. $A\in \F_{Bohr}$), if for any minimal (resp. minimal equicontinuous) system $(Y,S)$, there is $y\in Y$ such that $S^{n_i}y\ra y$ for some $\{n_i\}_{i=1}^\infty\subset A$. The Katznelson's question asks if any Bohr recurrence sequence is necessarily a Birkhoff one.

This problem remains unresolved and has attracted significant attention due to its deep connections with diverse mathematical areas. For related results and recent advances, we refer to \cite{Ka99, HY12, HSY16, HKM16,GKR22,DHM24, Gr23} and the references therein. Particularly, Katznelson \cite{Ka99} obtained a number-theoretical formulation of the question and related it with graph theory; Host, Kra and Maass \cite{HKM16} established that for a minimal pro-nilsystem the two types of recurrence sequences coincide, and if $A$ is a Birkhoff recurrence sequence for a minimal system $(Y,S)$, then $A$ is also a Birkhoff one for any minimal system $(Y',S)$ which is a proximal extension of $(Y,S)$; Glasscock, Koutsogiannis and Richter \cite{GKR22} demonstrated the equivalence of these sequences for specific skew product systems. In \cite{HY12} Huang and Ye formulated a group theoretical version of the question. Donoso, Heren\'andez and Maass \cite{DHM24} considered the question for $\mathbb{Z}^d$ actions and showed that Bohr recurrence sequences are recurrence sequences for $\mathbb{Z}^d$-Weyl systems.

Moreover, in \cite{HSY16} Huang, Shao and Ye introduced and investigated a high dimensional version of the question, and just recently Alweiss \cite{Al25} constructed a striking counterexample for the high dimensional version, which gives new evidence that Katznelson's question may have a negative answer.

\medskip
We may slightly modify the definition: Let $A\subset \N$ be called an {\it $R$-sequence (resp. $R_e$-sequence)}, denoted by $A\in \F_R$ (resp. $A\in \F_{R_e}$), if for any minimal (resp. minimal equicontinuous) system $(Y,S)$ there is $y\in Y$ such that $\{S^ny:n\in A\}$ is dense in $Y$.  In some personal communication Richter asked whether every $R_e$-sequence must necessarily be an $R$-sequence.

\medskip

Our subsequent discussion will reveal that Richter's question (was also listed as \cite[Question 17]{HSQY26}) is closely tied to the notion of {\it scattering},
introduced by Blanchard, Host, and Maass \cite{BHM}, as well as the notion of {\it weak scattering}, subsequently defined by Huang and Ye \cite{HY04} and systematically studied in \cite{HY12}. In this paper, we address Richter's question along with related problems in the context of ​​totally minimal systems​. Specifically, we will also consider subsets $A\subset \N$ satisfying the property that for every totally minimal system $(Y, S)$ and for every $y\in Y$, $\{S^ny: n\in A\}$ is dense in $Y$. Such a question is related to the notion of {\it disjointness} introduced by Furstenberg in \cite{F1}.

The motivation for extending Richter's question to totally minimal systems stems from the existence of some interesting sequences $A$ (see Corollary \ref{xuhui}) that satisfy the following property: for every totally minimal system $(Y,S)$, there exists a point $y\in Y$ whose orbit along $A$ is dense in $Y$, yet $A$ fails to be an $R$-sequence.

\subsection{Main results}
Our first main result is the characterization of $R$-sequences.

\medskip

\begin{thma}\label{thma}
 Let $A\subset \N$. The following statements are equivalent:
\begin{enumerate}
\item $A\in \F_R$.
\item $A$ is a shift-invariant Birkhoff recurrence sequence. i.e. for any $n\in\Z$, $A+n\in \F_{Bir}$.
\item For any syndetic set $S\subset \N$ and any $n\in\Z$, $A\cap (n+(S-S))\not=\emptyset.$
\item For any minimal system $(X,T)$ and any non-empty open set $U$ of $X$, $\cup_{a\in A}T^{-a}U$ is dense in $X$.
\end{enumerate}
\end{thma}

\medskip
Let $(X,T)$ and $(Y,S)$ be t.d.s. If $(X\times Y, T\times S)$ is transitive, then the set of transitive points \footnote{$x$ is a transitive point if
$\{T^nx:n\in\N\}$ is dense in $X$.} of $T\times S$ is a dense $G_\delta$ subset of
$X\times Y$. Then Veech's theorem (\cite[Proposition 3.1]{Veech70}) assets that there is a dense $G_\delta$ subset $\Omega\subset X$ such that for each $x\in \Omega$, there is $y\in Y$ such that $(x,y)$ is a transitive point. 
Remarkably, this set $\Omega$ can actually be chosen as the set of ​​transitive points of the system $(X,T)$ when
$(Y,S)$ is minimal.
\medskip

\begin{thma}\label{thmb}
 Assume that $(X,T)$ is transitive, $(Y,S)$ is minimal and $(X\times Y, T\times S)$  is transitive. Then for any transitive point $x\in X$ there is $y\in Y$ such that the orbit of $(x,y)$ is dense in $X\times Y$. 
 \end{thma}

\medskip
Theorem B has many corollaries.
Furstenberg \cite{F1} proved that if $(X,T)$ is minimal and weakly mixing then there is a dense $G_\delta$ subset $\Omega$ of $X$ such that for each $x\in \Omega$, $P[x]$  is a dense $G_\delta$ subset of $X$, where $P[x]$ is the set of all $y\in X$ which is proximal to $x$.
Akin and Kolyada \cite[Theorem 3.8]{AG03} improved the result by showing that if $(X,T)$ is weakly mixing then for each $x\in X$, $P[x]$ is dense in $X$ (later it was shown in \cite{HSY04} that $P[x]$ is a dense $G_\delta$ subset of $X$). The first corollary of Theorem \ref{thmb} is:

\medskip
\noindent {\bf Corollary B-1.} {\it Let $(X,T)$ be a minimal weakly mixing system. Then for any $x\in X$, the set $\{y\in X: (x,y)\ \text{has a dense orbit in}\ X\times X\}$ is a dense $G_\delta$ subset of $X$.}

\medskip
And a more general corollary is

\medskip
\noindent {\bf Corollary B-2.} {\it Let $(X,T)$ and $(Y,S)$ be minimal and $(X\times Y, T\times S)$ is transitive. Then for each $x\in X$, $\{y\in Y: (x,y)\ \text{has a dense orbit in}\ X\times Y\}$ is a dense $G_\delta$ subset of $Y$.}

\medskip

For a transitive system $(X,T)$, we say a transitive point $x$ is a {\it scattering point} if  for any minimal system $(Y,S)$ there is $y\in Y$ such that $(x,y)$ has a dense orbit. It can be shown that (Theorem \ref{equi-s-point}): if $(X,T)$ is transitive, then a transitive point $x$ is a scattering point if and only if for any open non-empty set $U$ of $X$, $N(x,U)\in \F_R$, where $N(x,U)=\{n\in\N:T^nx\in U\}$.

\medskip
It was shown in \cite{HY02} that $(X,T)$ is scattering if and only if for any open non-empty sets $U$ and $V$ of $X$, $N(U,V)$ is a Birkhoff recurrence sequence, which was used to construct an explicit system that is scattering and not weakly mixing, where $N(U,V)=:\{n\in\N: U\cap T^{-n}V\not=\emptyset\}$. Below, we establish additional equivalent formulations of scattering.

\medskip

\begin{thma}\label{thmc}
 Let $(X,T)$ be a t.d.s. The following statements are equivalent
\begin{enumerate}
\item $(X,T)$ is scattering.

\item For any transitive point $x\in X$ and for any open non-empty subset $U$ of $X$, $N(x,U)\in \F_R$.

\item For any transitive point $x\in X$ and any minimal system $(Y,S)$, there is $y\in Y$ such that the orbit of $(x,y)$ is dense in $X\times Y$.
\end{enumerate}
\end{thma}

\medskip
The above theorem gives a lot of examples of $R$-sequences. For instance, if $(X,T)$ is weakly mixing then $N(x,U)\in \F_R$ for any transitive point $x$ and any non-empty open subset $U$ of $X$, since weak mixing implies scattering.

\medskip

It is easy to show that $(X,T)$ is weakly scattering if and only if for any open non-empty set $U$ of $X$, $N(x,U)\in \F_{R_e}$ (Theorem \ref{F-ER}). So combining this result with Theorem \ref{thmb} and earlier work of Huang and Ye \cite{HY12} we deduce that

\medskip
\noindent {\bf Corollary C-1.} {\it If scattering and weak scattering are distinct properties, then both Richter's question and Katznelson's question have a negative answer.}

\medskip
Let $(X,T)$ and $(Y,S)$ be t.d.s. We say that $J\subset X\times Y$ is a {\it joining} if $J$ is non-empty, closed and invariant ($T\times S(J)\subset J$); and $(X, T)$ is disjoint from $(Y,S)$ (denoted by $X\perp Y$) if each joining is equal to $X\times Y$. It is easy to see that if $X\perp Y$ then one of them must be minimal (\cite[Theorem II.1]{F2}).

One of the fundamental questions related to disjointness, posed by Furstenberg \cite{F2}, is to characterize all t.d.s. that are disjoint from ​​every minimal system​​. This question has been investigated by Huang and Ye in \cite{HY05}, as well as by Huang, Shao, and Ye in \cite{HSY20}, where some non-intrinsic characterizations of systems disjoint from all minimal systems were presented. For other related works, see \cite{LYY15, Op19}. Notably, in a recent preprint, Huang, Shao, Xu, and Ye in \cite{HSXY} obtained an ​​intrinsic characterization​​ of systems that are disjoint from all minimal systems, answering the question completely. We remark that Furstenberg's question can be extended for group actions, see \cite{XY22, GTWZ21} for the progress on this respect.

\medskip

Let $\mathcal{M}$ (resp. $\mathcal{M}_{TM}$) be the collection of all minimal (resp. totally minimal) systems and $\mathcal{M}^{\perp}$ (resp. $\mathcal{M}_{TM}^{\perp}$) be the collection of systems that disjoint from all minimal  (resp. totally minimal) systems . A transitive system $(X,T)$ is an $M$-system if the set of minimal points is dense in $X$. It was proved in \cite[Remark 4.8]{DSY12} that if $(X,T)$ is  transitive and $X\in \mathcal{M}_{TM}^\perp$ then $(X,T)$ is an $M$-system. 

For a minimal system $(X,T)$, each minimal subset in $(2^X,T)$ is called a {\it quasi-factor}, where $2^X$ is the collection of all non-empty closed subsets of $X$ equipped with the Hausdorff metric.

Let $(X,T)$ be an $M$-system. We denote  the closure of $\bigcup\{M\subset X: M \text{ is minimal and } M\perp \mathcal{M}_{TM}\}$ in  $X$ 
by $X_{TM}^\perp$. Note that in \cite[Theorem 6.3]{HSXY} the authors proved that a transitive system $(X,T)\in \mathcal{M}^\perp$ if and only if $X$ is an $M$-system and it is disjoint from any minimal subset of $X$. Since there are minimal t.d.s. disjoint from all totally minimal systems (Theorem \ref{minimal-disjoint}),
the formulation of the following theorem is more involved than the setting considered in \cite[Theorem 6.3]{HSXY}.

\medskip

\begin{thma}\label{thmd}
 Let $(X,T)$ be a transitive t.d.s. Then the following statements are equivalent
\begin{enumerate}
\item $(X,T)\perp \mathcal{M}_{TM}.$
\item Either there is a sequence of  minimal sets $\{M_i\}_{i=1}^\infty\subset X$ with $\overline{\bigcup_{i=1}^{\infty}M_i}=X$ such that $M_i\perp \mathcal{M}_{TM}$ for each $i\in \N$, or $X_{TM}^\perp$ is nowhere dense and there is a sequence of minimal sets $\{N_i\}_{i=1}^\infty\subset X$ with $\overline{\bigcup_{i=1}^{\infty}N_i}=X$ such that $N_i\not\perp\mathcal{M}_{TM}$ and $X$ is disjoint from all totally minimal quasifactors of $N_i$ for each $i\in\mathbb{N}$. 
\end{enumerate}
\end{thma}

\medskip
\begin{rem}
It is relatively easy to construct an $R$-sequence with zero Banach density \footnote{The upper Banach density $A\subset \N$ is defined by
$\limsup_{\{I_i\}}\frac{1}{|I_i|} |A\cap I_i|$, where $\{I_i\}$ runs over intervals of $\N$.}. If $p$ is an integral polynomial vanish at zero, then $\{p(n):n\in\N\}$ is an $R$-sequence for all totally minimal systems, see \cite{GHSWY20, QJ23}. At the same time, the sequence $A\subset \N$ satisfying that for any minimal
(resp. totally minimal) system $(Y,S)$ and any $y\in Y$, $\{S^ny: \ n\in A\}$ is dense must have positive upper Banach density, and in fact $A$ is piecewise syndetic (meaning $A$ is ``large"), see  \cite[Theorem 2.4]{HY05} by Huang and Ye, and \cite[Proposition 4.4]{DSY12} by Dong, Shao and Ye.
In a recent preprint, Glasscock and Le \cite{GL24} obtained a completely characterization of the  sequence $A$ satisfying that for any minimal
system $(Y,S)$ and any $y\in Y$, $\{S^ny:n\in A\}$ is dense.
\end{rem}

\begin{rem} If we consider the sequence $A\subset \N$ such that for each minimal system $(Y,S)$ and each $y\in Y$ there is a sequence $\{n_i\}_{i=1}^\infty\subset A$ with $S^{n_i}y\ra y$, then this property is closely related to the notion of {\it weak product recurrence}.
Recall that  $x$ is a weak product recurrent point if $(x,y)$ is recurrent for any minimal point $y$.

The concept of {\it ​​product recurrence}​​ was originally introduced by Auslander and Furstenberg \cite{AF94}, while the refined notion of weak product recurrence was later defined by Haddad and Ott \cite{HO08} who showed that there are weak product recurrent points that are not minimal. Moreover, Glasner and Weiss \cite{GW17} constructed a minimal (not distal) weak product recurrent point (in fact it is the transitive point of the so-called double minimal systems). Dong, Shao and Ye \cite[Proposition 4.4]{DSY12} demonstrated that any such sequence $A$ must necessarily be ​​piecewise syndetic​​. Furthermore, Glasscock and Le \cite{GL24} obtained a complete characterization of these sequences, providing a comprehensive understanding of their structural properties.
\end{rem}

\begin{rem} In the current paper, we characterize transitive systems that are disjoint from all totally minimal systems. In a forthcoming paper by the same authors, we will characterize topological dynamical systems  disjoint from all totally minimal systems, and measure-preserving systems disjoint from all totally ergodic systems. We note that measure-preserving systems disjoint from all ergodic systems have been characterized in \cite{GLR24, HSXY} and for a short proof, see \cite{GW24}.
\end{rem}

\subsection{The organization of the paper}
In Section 2, we provide the necessary preliminaries following the introduction. Section 3 is devoted to a brief discussion of Katznelson's question, the related one for totally minimal systems and the proof of Theorem \ref{thma} as a worm up. In Section 4, we investigate $R$-sequences and establish the proofs of Theorems \ref{thmb} and \ref{thmc}. Section 5 focuses on systems whose product with all totally minimal systems are transitive, and Section 6 contains the proof of Theorem \ref{thmd}, that is we characterize transitive systems disjoint from all totally minimal systems. Finally, summary of the results is given in Section 7.

\bigskip
\noindent {\bf Acknowledgement:} We would like to thank
W. Huang and S. Shao for useful discussions. Thank J. Li for the careful reading which improves the writing of the paper, and J. Griesmer for providing us many useful references.

\section{Preliminary}

In this section we give some necessary notions and some known facts which
we will use later.

\subsection{The general notions}
A {\it topological dynamical system} (t.d.s. or system for short) is a compact metric space $X$ together with a continuous onto map $T:X\to X$. For open sets $U,V\subset X$, let
\[ N(U,V)=\{n\in\N:U\cap T^{-n}V\not=\emptyset\}.\]

For $x\in X$ and an open set $U\subset X$, let
\[N(x, U)=\{n\in\N: \ T^{n}x\in U\}.\]
A system $(X, T)$ is {\it transitive} if for  non-empty  open subsets $U$ and $V$ of $X$, $N(U,V)\not=\emptyset$;  is {\it weakly mixing} if
$(X\times X,T\times T)$ is transitive.

For a system $(X,T)$, the {\it orbit of $x$} is the set $\{T^nx:x\in\N\}$, denoted by $orb(x,T)$, or $\O(x,T)$ whose closure is denoted by  $\overline{\O}(x,T)$. It is not hard to show that if $(X,T)$ is transitive then the set of transitive points, i.e. the points with dense orbits, forms a dense $G_\delta$ set of $X$. We use $Tran_{T}(X)$ (or $Tran_{T}$) to denote the set of transitive points of $(X,T)$.

\medskip
A system $(X,T)$ is {\it equicontinuous} if for each $\ep>0$ there is $\delta>0$ such that $\rho(x,y)<\delta$ implies that
$\rho(T^nx,T^ny)<\ep$ for each $n\in\N$. We say that $x,y\in X$ are {\it proximal} if $\inf\limits_{n\in\N}\rho(T^nx,T^ny)=0$, otherwise we say that $x,y$
are {\it distal}. Moreover, we say that $(X,T)$ is {\it distal} if for any $x\neq y$ in $X$ then $x,y$ are distal; is {\it proximal} if for any $x,y\in X$, $x,y$ are proximal.

\subsection{Minimality}
\subsubsection{General notions}
We say that $(X,T)$ is {\it minimal} if the orbit of each $x\in X$ is dense. It is a fundamental fact that each t.d.s. $(X,T)$ has a subsystem which is minimal.
We say that $x$ is a {\it minimal point} in a t.d.s. $(X,T)$ if the orbit closure of $x$ is a minimal subset of $X$. A minimal system $(X,T)$ is {\it totally minimal} if $(X,T^n)$ is minimal for any $n\in\N$.

It is known if $(X,T)$ is minimal then $T$ is necessarily surjective but may not be injective. However, we have
\begin{lem}\cite{KTS01}\label{min map}
If $(X,T)$ is minimal then
\begin{enumerate}
\item[(1)] $T$ is semi-open, that is the image of any nonempty open set under $T$ has nonempty interior.
\item[(2)] The set $\{x\in X: |T^{-1}\{x\}|=1\}$ is a dense $G_{\delta}$ subset of $X$.
\end{enumerate}
\end{lem}

\begin{lem} Let $(X,T)$ be minimal. Then $(X,T)$ is totally minimal if and only if for any $n\in\N$ and any $x\in X$ there is a sequence $\{k_i\}_{i=1}^\infty\subset \N$ such that $T^{nk_i+1}x\ra x$.
\end{lem}
\begin{proof} The necessity is  clear.  Suppose that for any $n\in\N$ and any $x\in X$ there is a sequence $\{k_i\}_{i=1}^\infty\subset \N$ such that $T^{nk_i+1}x\ra x$. Then for each $x\in X$, $T^{n-1}x\in \overline{\O}(x, T^{n})$. Inductively, we have that
\[T^{n-i}x\in \overline{\O}(T^{n-(i-1)}x, T^{n})\subset \overline{\O}(T^{n-(i-1)}x, T^{n})\subset\cdots\subset \overline{\O}(x, T^{n}),\ \forall i=1,2,\ldots,n.\]
For any $x,y\in X$, $y\in \overline{\O}(x,T)=\cup_{i=0}^{n-1}\overline{\O}(T^{i}x, T^n)$. Thus $y\in \overline{\O}(x, T^n)$. This shows that $(X,T)$ is totally minimal.
\end{proof}

If $(X,T)$ is a minimal system then we use $(X_{eq}, T)$ to denote its maximal equicontinuous factor, which is deduced by the regionally proximal relation.

The collection of all minimal (resp. minimal equicontinuous and minimal distal) systems is denoted by $\mathcal{M}$ (resp. $\mathcal{M}_E$
and $\mathcal{M}_D$). The collection of all totally minimal (resp. totally minimal equicontinuous and totally minimal distal) systems is denoted by $\mathcal{M}_{TM}$ (resp. $\mathcal{M}_{TE}$ and $\mathcal{M}_{TD}$).

A transitive system is
\begin{itemize}
 \item an $M$-system if the set of minimal points is dense;
 \item a $P$-system if the set of periodic points is dense;
 \item an $E$-system if it admits an invariant probability measure with full support.
\end{itemize}
It is clear that an $M$-system must be an $E$-system. The following is a well-known result see
for instance \cite[Lemma 2.1]{HY05}.

\begin{lem} \label{char of M-sys}
 A transitive t.d.s. $(X,T)$ is an $M$-system if and only if for each transitive point $x$ and each neighborhood $U$ of $x$, $N(x,U)$
is piecewise syndetic.
\end{lem}
See Subsection \ref{sb2.5} for the definition of piecewise syndeticity.

\subsubsection{PI systems}
Let $\pi:X\ra Y$ be an extension. We say that {\it $\pi$ is equicontinuous} if for any  any $\ep>$ there is $\delta>0$ such that $\rho(x,y)<\delta$
and $\pi(x)=\pi(y)$ implies that $\rho(T^nx,T^ny)<\ep$ for any $n\in \N$; and {\it $\pi$ is proximal} if $\pi(x)=\pi(y)$ implies that $x,y$ are proximal.

{\it A  minimal system is PI} if it is a factor of an inverse limit of equicontinuous or proximal extensions starting from the trivial system.

\subsection{Weak disjointness and disjointness}
Let $(X,T)$ and $(Y,S)$ be two systems. We say that they are {\it weakly disjoint}, denoted by $X\curlywedge Y$,  if the product system $(X\times Y, T\times S)$ is transitive. A {\it joining} of $X$ and $Y$ is a closed subset $J\subset X\times Y$ such that $T\times S(J)\subset J$ and the projections to $X$ and $Y$ are full. If $X\times Y$ is the unique joining of $X$ and $Y$ then we say that they are {\it disjoint} and we denote it by $X\perp Y$.

We use $\mathcal{M}^{\curlywedge}$ (resp. $\mathcal{M}_{E}^{\curlywedge}$) to denote the collection of systems that are weakly disjoint from all minimal (resp. minimal equicontinuous) systems. A t.d.s. is scattering if it is in $\mathcal{M}^\curlywedge$,
it is weakly scattering if it is in $\mathcal{M}_E^\curlywedge$.   

It was shown by Blanchard, Host and Maass \cite{BHM} that a scattering system is in $\mathcal{M}_D^\perp$, and  by Huang and Ye \cite{HY12} that a weakly scattering system is also in $\mathcal{M}_D^\perp$. Considering totally minimal distal systems, we also have


\begin{prop}\label{simpleproof} If $(X,T)\in \mathcal{M}_{E}^\curlywedge$ (resp. $\mathcal{M}_{TE}^\curlywedge$), then it is in $\mathcal{M}_{D}^\perp$ (resp. $\mathcal{M}_{TD}^\perp$).
\end{prop}

\begin{proof}[Sketch of the proof]
We only need to show that $(X,T)\in \mathcal{M}_{TE}^\curlywedge$ implies $(X,T)\in \mathcal{M}_{TD}^\perp$. First, it is clear that $\mathcal{M}_{TE}^\curlywedge=\mathcal{M}_{TE}^\perp$. Let $(X,T)\in \mathcal{M}_{TE}^\curlywedge$ and let $(Y,S)\in \mathcal{M}_{TD}$. Since the disjointness is preserved under inverse limits, it only needs to show that if $\pi: Y\rightarrow Z$ is an isometric extension and $X\perp Z$ then $X\perp Y$. By Lemma \ref{weaklydsijointmm}, one has $X\curlywedge Y$ since $X\perp Y_{eq}$. Now for any joining $J$ of $X$ and $Y$, $(\id\times \pi)(J)=X\times Z$. Since $\id\times \pi: X\times Y\rightarrow X\times Z$ is a group extension and $X\times Y$ is transitive, it follows from \cite[Lemma 4.4]{BHM} that $J=X\times Y$ and hence $X\perp Y$.
\end{proof}

\begin{lem}\cite[Proposition 4.14]{HY12}\label{E-w-s}
If $(X, T)$ is an $E$-system, then $(X, T)$ is weakly scattering if and only if  it is weakly mixing.
\end{lem}

\subsection{Almost equicontinuity}
The notion of {\it almost equicontinuity} was introduced by Glasner and Weiss \cite{GW93}.

Let $(X,T)$ be a t.d.s. $x\in X$ is an {\it equicontinuous point} if for any $\ep>0$ there is $\delta>0$ such that $\rho(x,y)<\delta$ implies that
$\rho(T^nx,T^ny)<\ep$ for any $n\in\N$ and $y\in X$. $(X,T)$ is {\it almost equicontinuous} if there is at least one equicontinuous point.

We remark that
if $(X,T)$ is transitive and almost equicontinuous then each transitive point of $X$ is an equicontinuous point.

\begin{lem}\cite[Theorem 4.2]{AG01}\label{almost-eqi} Let $(X,T)$ be a transitive t.d.s. and $x$ be a transitive point. Then $(X,T)$ is almost equicontinuous if and only if for each $\ep>0$ there is $\delta>0$ such that $\rho(T^ix,x)<\delta$ implies that $\rho(T^{i+j}x,T^jx)<\ep$ for any $j\in \N$.
\end{lem}

\begin{lem}\label{AG-exa}\cite[Corollary 4.14]{AG01} or \cite[Remark 5.5]{HY12} There is an almost equicontinuous system (transitive but not minimal) which is scattering.
\end{lem}

Since an almost equicontinuous system is not weakly mixing, this shows that scattering and weak mixing are different properties.

\begin{lem}\label{GW93}\cite[Theorem 4.6]{HY02}
If $(X,T)$ is a transitive almost equicontinuous system, not minimal then the set of minimal points is not dense in $X$.
\end{lem}

Thus, we have

\begin{cor}\label{dis-equi-dis} Disjointness from all minimal (resp. totally minimal) systems and disjointness from all minimal equicontinuous (resp. totally minimal equicontinuous) systems are different properties.
\end{cor}
\begin{proof} Let $(X,T)$ be the almost equicontinuous system  described in Lemma \ref{AG-exa}. This system is disjoint from all minimal equicontinuous systems, as established in Theorem \ref{simpleproof}. However, it cannot be disjoint from all minimal systems, since such a property would imply that $(X,T)$ is an $M$-system (see \cite[Theorem 2.6]{HY05}), a conclusion that contradicts the result of Lemma \ref{GW93}.

To show the totally minimal case, we only need to use the arguments in \cite[Theorem 2.6]{HY05} with \cite[Remark 4.8]{DSY12}.

\end{proof}


\subsection{Generic eigenvalues} For a transitive system $(X,T)$ let $Eig(X,T)$  be the set of all generic eigenvalues of $(X,T)$.
Recall that if $(X, T)$ is a transitive t.d.s., then a complex number $\lambda\in \mathbb{C}$ called a {\it generic eigenvalue} of $(X, T)$ if there exists  a continuous nonzero function $f$ on $Tran_T (X)$ with $T\circ f = \lambda f$. The function $f$ is called the generic eigenfunction (associated with $\lambda$). 

The following lemma was proved in \cite[Theorem 2.7]{HY12}.

\begin{lem}\label{HY-scattering} Let $(X, T)$ be a transitive t.d.s. Then the following statements are equivalent:
\begin{enumerate}
\item $(X, T)$ is weakly scattering.
\item $(X, T)$ is totally transitive and weakly disjoint from any irrational rotation $(\T, R_\lambda)$.
\item $(X, T)$ has no non-trivial generic eigenvalues, i.e. $Eig(T) = \{1\}$.
\end{enumerate}
\end{lem}

The following lemma is a part of \cite[Proposition 3.7]{HY12}.

\begin{lem}\label{eign-genreic} Let $(X,T)$ be transitive and $(Y,S)$ be minimal and equicontinuous. Then $X\perp Y$ if and only if $Eig(X,T)\cap Eig(Y,S)=\{1\}$.
\end{lem}

Similar to Lemma \ref{HY-scattering} we have
\begin{thm}\label{HY-ext}Let $(X, T)$ be a transitive t.d.s. Then the following statements are equivalent:
\begin{enumerate}
\item $(X, T)$ is weakly disjoint from all totally minimal equicontinuous systems.
\item $(X, T)$ is weakly disjoint from any irrational rotation $(\T, R_\lambda)$.
\item Each generic eigenvalue of $(X, T)$ is a root of the unit.
\end{enumerate}
\end{thm}
\begin{proof} It follows by Lemma \ref{eign-genreic}.
\end{proof}

\subsection{Some subsets of $\N$}\label{sb2.5}\
\medskip

A subset $S$ of $\mathbb{N}$ is {\it syndetic} if it has a bounded
gap, i.e., there is $N \in \mathbb{N}$ such that $\{i, i+1, \ldots,
i+N\} \cap S \neq \emptyset$ for every $i \in \mathbb{N}$. A subset $S\subseteq \N$
is {\it thick} if it contains arbitrarily long runs of integers, i.e., for every $n \in \mathbb{N}$ there exists some $a_n
\in \mathbb{N}$ such that $\{a_n, a_n+1, \ldots, a_n+n\} \subseteq S$.

A subset $S$ of $\mathbb{N}$ is {\it piecewise syndetic} if it is the intersection of a syndetic set with a thick set;
and it is {\it thickly syndetic} if for each $n\in\N$ there is a syndetic
subset $\{w^n_1,w^n_2, \ldots\}$ of $S$ such that
$\{w^n_i,w^n_i+1,\ldots,w^n_i+n\}\subseteq S$ for each $i\in\N$.
Denote by $\F_{ps}$ (resp. $\F_{ts}$) the family of all piecewise syndetic set (resp. all thickly syndetic sets). It is clear that if $F_1,F_2\in \F_{ts}$ so is $F_1\cap F_2$. That is, $\F_{ts}$ is a filter.

Let $\F$ be a family of $\N$, then the family $\{B\subset \N: B\cap A\not=\emptyset\ \text{for each}\ A\in \F\}$ is denoted by $\F^*$.

Let $F=\{n_i\}_{i=1}^d$ be a subset of $\N$, then $\{n_{i_1}+\cdots+n_{i_j}: i_1, \ldots, i_j\in F, j\le d\}$ is called a finite $IP$-set (generated by $F$), and
if $F=\{n_i\}_{i=1}^\infty$, then $\{n_{i_1}+\cdots+n_{i_j}: i_1, \ldots, i_j\in F, j\in\N\}$ is called an $IP$-set (generated by $F$), denoted by $FS(F)$.

\section{Recurrence  and weak product recurrence}

In this section, we recall some known definitions of some kinds of recurrence sequences and give  the definition of $R$-sequences.

\subsection{Recurrence sequence}

Recall that $A\subset \N$ is a {\it topological recurrence sequence} is if for any t.d.s. $(X,T)$, there is some $x\in X$ such that for any neighborhood $U$ of $x$, there is some $n\in A$ such that $T^nx\in U$. Further,  $A\subset \N$ is a {\it Birkhoff recurrence sequence} 
if for any minimal system $(X, T)$, there is some $x\in X$ such that  for each neighborhood $U$ of $x$, there is $n\in A$ such that $T^nx\in U$. We use $\mathcal{F}_{Bir}$ to denote the collection of Birkhoff recurrence sequences.

Let $(X,T)$ be a system and $A\subset \N$.  A point $x\in X$ is {\it recurrent along $A$} is for any neighborhood $U$ of $x$, there is some $n\in A$ such that $T^nx\in U$. This is equivalent that there is an increasing subsequence $(n_i)$ in $A$ such that $T^{n_i}x\rightarrow x$ as $i$ tends to $\infty$. We say $x$ is recurrent if it is recurrent along $\N$.

The following equivalent characterizations is classical, see for instance \cite{HKM16,HSY16}. 
\begin{thm}\label{equi of rec}
Let $A\subset \N$. Then the following statements are equivalent.
\begin{enumerate}
\item $A$ is a topological recurrence set.
\item $A$ is a Birkhoff recurrence sequence.

\item For any minimal system $(X, T)$, there is a dense $G_\delta$ subset $\Omega$ such that for each
$x\in \Omega$ and each neighborhood $U$ of $x$, there is $n\in A$ such that $T^nx\in U$.

\item For each minimal system $(X,T)$ and any non-empty open set $U$, $A\cap N(U,U)\not=\emptyset.$

\item For any syndetic set $S$ of $\N$, $A\cap (S-S)\not=\emptyset.$
\end{enumerate}
\end{thm}

\subsection{$R$-sequences and the proof of Theorem \ref{thma}} Recall that a subset $A\subset \mathbb{N}$ is an $R$-{\it sequence} if for any any minimal system $(X,T)$, there is some $x\in X$ such that the orbit of $x$ along $A$ is dense in $X$, that is $orb_{A}(x):=\{T^{n}x:\ n\in A\}$ is dense.
The collection of all $R$ sequences is denoted by $\F_{R}$.

The following equivalent characterizations  is clear. For the completeness, we give a proof.
\begin{lem}\label{s=dence}
Let $A\subset \N$ and $(X,T)$ be minimal. Then the following statements are equivalent.
\begin{enumerate}
 \item There is $x\in X$ such that $\{T^nx:n\in A\}$ is dense in $X$.
 \item For each open non-empty subset $V$, $\cup_{n\in A}T^{-n}V$ is dense in $X$.
 \item There is a dense $G_\delta$ subset $\Omega$ such that for each $y\in \Omega$, $\{T^ny:n\in A\}$ is dense in $X$.
 \end{enumerate}
\end{lem}
\begin{proof} (1) $\Rightarrow$ (2) First note that for any $n\in \N$, $T^nx$ also has dense orbit along $A$. This implies that for any open non-empty subsets $U,V$ we have that $U\cap (\cup_{n\in A}T^{-n}V)\not=\emptyset$ which in turn implies that $\cup_{n\in A}T^{-n}V$ is dense in $X$.

(2) $\Rightarrow$ (3)  For a base $\{V_i\}$ of the topology of $X$, set $\Omega=\cap_{i=1}^\infty (\cup_{n\in A}T^{-n}V_i)$. Then $\Omega$ is the set we need. It is clear that (3) implies (1).
\end{proof}

\medskip

\medskip

Let $\F_M$ be the family generated by $N(U,V)$, where $(X,T)$ is a minimal system $(X,T)$ and $U,V$ are open non-empty sets of $X$. Theorem \ref{thma} is contained in the following theorem.

\begin{thm}\label{rdop-char}
For $A\subset \mathbb{N}$, the following statements are equivalent
\begin{enumerate}
\item $A\in \F_R$.
\item $A\in \F_M^*$, i.e. for any minimal system $(X,T)$ and open non-empty sets $U$ and $V$, $A\cap N(U,V)\not=\emptyset.$
\item For any minimal system $(X,T)$ and any non-empty open subsets $U$, $\cup_{a\in A}T^{-a}U$ is dense in $X$.
\item For any minimal system $(X,T)$, any $n\in\Z$ and any non-empty open subset $U$, $(A+n)\cap N(U,U)\not=\emptyset$.
\item $A$ is a shift-invariant Birkhoff recurrence sequence. i.e. for any $n\in\Z$, $(A+n)\cap \N$ is a Birkhoff recurrence sequence.
\item For any syndetic set $S\subset \N$ and any $n\in\Z$, $(A+n)\cap (S-S)\not=\emptyset.$
\end{enumerate}
\end{thm}
\begin{proof} It follows from Lemma \ref{s=dence}  and  Theorem \ref{equi of rec} that (1) $\Leftrightarrow$ (3) and (5) $\Leftrightarrow$ (6), respectively. Next we show that (1) $\Rightarrow$ (2) $\Rightarrow$ (4) $\Rightarrow$ (5) $\Rightarrow$ (1).

(1) $\Rightarrow$ (2) Assume that $A\in \F_R$ and $(X,T)$ is a minimal system. Then there is $x\in X$ such that $\{T^nx:n\in A\}$ is dense in $X$.
In fact, by Lemma \ref{s=dence}, such $x$ forms a dense $G_\delta$ set. Now let $U,V$ be open non-empty sets of $X$. Then there is $y\in U$
such that $\{T^ny:n\in A\}$ is dense in $X$. It implies that there is $n\in A$ such that $T^ny\in V$ which in turn implies $U\cap T^{-n}V\not=\emptyset$ i.e. $A\cap N(U,V)\not=\emptyset$ or $A\in \F_M^*$.

\medskip



(2) $\Rightarrow$ (4) Let $U$ be a nonempty open set of $X$. If $n\geq 0$ then $T^{-n}U$ is also non-empty open set. It follows from that  (2) that there is some $a\in A$ such that $T^{-a}T^{-n}U\cap U\neq \emptyset$. This implies that $a+n\in N(U,U)$ and hence $(A+n)\cap N(U,U)\neq\emptyset$. If $n<0$ then it follows from Lemma  \ref{min map} that  the interior of $T^{-n}U$, which is denoted by $V$, is non-empty. By (2), there is some $a\in A$ such that $T^{-a}V\cap U\neq \emptyset$. Thus $T^{-a}T^{-n}U\cap U\neq\emptyset$ and we also have $(A+n)\cap N(U,U)\neq\emptyset$.

\medskip
(4) $\Rightarrow$ (5) Let $(X,T)$ be a minimal system and  $n\in \Z$. Now for any $\epsilon>0$, we can take a non-empty open set $U$ with $\diam(U)<\epsilon$. Then there is some $k\in (A+n)\cap \N$ such that $T^{-k}U\cap U\neq \emptyset$.  Thus there is some $x\in U$ such that $\rho(T^{k}x,x)<\epsilon$. In particular,
\[ \Omega_{\epsilon}:=\{x\in X: \exists k\in (A+n)\cap \N \text{ such that } \rho(T^{k}x,x)<\epsilon\}\neq\emptyset, \ \forall \epsilon>0.\]
Further, by the minimality of $(X,T)$, it is easy to see that $\Omega_{\epsilon}$ is open dense.  By the Baire's Theorem, the set $\Omega=\bigcap_{i=1}^{\infty}\Omega_{1/i}$ is nonempty.  It is clear that any point in $\Omega$ is recurrent along $(A+n)\cap \N$.

\medskip
(5) $\Rightarrow$ (1) Let $(X,T)$ be a minimal system. By Theorem \ref{equi of rec} , the set $\Omega_n$ consisting of recurrent points along $A+n$ is a dense $G_{\delta}$ set for each $n\in \N$. Thus $\Omega=\bigcap_{m=1}^{\infty}T^{-m}(\bigcap_{n=1}^{\infty})$ is also a dense  $G_{\delta}$ set. Now take any $x\in \Omega$. For any nonempty open set $U$ in $X$, it follows from the minimality that there is some $n\in\N$ such that $T^{n}x\in U$.  Note that $T^{n}x\in \Omega$ and then $T^{n}x$ is recurrent along $(A-n)\cap \N$. In particular, there is some $a\in A$ such that $T^{a}x=T^{a-n} T^{n}x\in U$. Thus $\{T^{a}x: a\in A\}$ is dense in $X$ as $U$ is chosen arbitrarily.
\end{proof}



\subsection{Product recurrence}
By \cite[Proposition 2.14]{F1} it is easy to see that if $(X,T)$ is transitive and $x\in Tran_T$, then for any minimal system $(Y,S)$, $R_x=\{y\in Y: (x,y)\ \text{is recurrent}\}$
is a dense $G_\delta$ set of $Y$. Particularly, when $(X,T)$ is minimal then $R_x$ is a such set for each $x\in X$.

\medskip

We say that a transitive $(X,T)$ has {\it weak product property} if there is $x\in Tran_T$ such that for each minimal system $(Y,S)$ and any $y\in Y$, $(x,y)$ is recurrent. Such a point $x$ is called a {\it weakly product recurrent} point. Note that it was shown in \cite{DSY12} that a transitive system with weak product property must be an $M$-system.

Recall that we say that $A\subset \N$ is a {\it pointwise recurrence sequence} if for any minimal system $(Y,S)$ and any $y\in Y$, there is $B=\{b_i\}\subset A$ such that $S^{b_i}y\ra y$.

The connection of the above two notions is

\begin{lem}
Let $(X,T)$ be a transitive system. Then $x\in Tran_T$ is a weakly product recurrent point if and only if for any neighborhood $U$ of $x$, $N(x,U)$ is a pointwise recurrence sequence.
\end{lem}
\begin{proof} First assume that $x$ is a weakly product recurrent point, and $U$ is a neighborhood of $x$. Let $(Y,S)$ be minimal and $y\in Y$.
Let $V_i$ be a sequence of neighborhood of $y$ with $\diam(V_i)\ra 0$. Then for any $i\in \N$ there is $n_i\in N(x,U)$ such that $S^{n_i}y\in V_i$
which implies that $S^{n_i}y\ra y$. That is, $N(x,U)$ is a pointwise recurrence sequence.

Conversely, assume that for any neighborhood $U$ of $x$, $N(x,U)$ is a pointwise recurrence sequence. Then it is clear that for any minimal point $y$, $(x,y)$ is a recurrence point. That is $x$ is a weakly product recurrent point.
\end{proof}

\begin{prop}\label{wprd-equ} There is $A\subset \N$ such that $A$ is a pointwise recurrence sequence for all minimal equicontinuous systems, but $A$ fails to be
a pointwise one for all minimal systems.
\end{prop}
\begin{proof} 


We take an IP set $A$ that is not piecewise syndetic. Such set indeed exists as shown in \cite{BR02} that there  is some translate of square-free set is an IP set but not piecewise syndetic. Then $B:=\mathbb{N}\setminus  A$ is thickly syndetic.  Now it follows from \cite[Theorem 2.4]{HY05} that there is some nonempty subset $C\subset B$ such that ${\bf 1}_{C}$ is a minimal point in $\Sigma_2:=\{0,1\}^{\mathbb{Z}_{+}}$. But then $N({\bf 1}_{C}, U)=C\subset B$, where $U=\{(x_n)\in\Sigma_2:\ x_0=1\}$. Thus $A$ is not pointwise recurrent for a minimal point ${\bf 1}_{C}$. On the other hand, $A$ is a pointwise recurrent for any distal minimal point and hence for minimal equicontinuous system.
\end{proof}

We remark that  the system $Y$ can be minimal weakly mixing and hence totally minimal. Thus, there is $A\subset \N$ such that $A$ is a pointwise recurrence sequence for all minimal equicontinuous systems, but $A$ fails to be
a pointwise one for all totally minimal systems.

\medskip

The following theorem gives some ways to construct pointwise recurrence sequence.
\begin{thm}\label{wpp with dense distal}
\begin{enumerate}
\item[(1)] Let $(X,T)$ be an $M$-system (not minimal with dense minimal sets $\{M_i\}$) and $x\in Tran_T$. Let $U_i$ be a sequence of  neighborhoods of $x$ with $diam(U_i)\ra 0$. If $x$ is proximal to any point $z\in U_i\cap M_i$ if it is not empty, then for each neighborhood $U$
of $x$, $N(x,U)$ is a pointwise recurrence sequence.


\item[(2)] Let  $(X,T)$ be transitive with dense set of distal subsets $(M_i)_{i\in\N}$ and $x\in Tran_T$. If
there is $p_j\in \cup_{i\in\N} M_i$ with $p_j\ra x$ and each $(x,p_j)$ is proximal, then $x$ is a weakly product recurrent point.
\end{enumerate}
\end{thm}

\begin{proof}
(1) Let $(Y,S)$ be minimal and $y\in Y$. Assume that $U_i\cap M_i\not=\emptyset$ for some $i\in \N$.
Then $M_y=\{z\in M_i: (z,y)\ \text{is minimal}\}$ is dense in $M_i$. Take $z_i\in M_y\cap U_i$ and assume that $B_{2\ep}(z_i)\subset U_i$.
It is clear that
$A=\{n\in\N:\rho(T^nx,T^nz_i)<\ep\}\in \F_t$ and $B=\{n\in\N: \rho(T^nz_i,z_i),\rho(S^ny,y)<\ep\}\in \F_s$. Since $A\cap B\not=\emptyset$,
and $x$ is a limit of minimal points from $\cup_{i\in\N}M_i$, the conclusion follows.

\medskip
(2)  Suppose that there is a sequence $(p_j)$ of distal points converging to $x$ such that each $(x,p_j)$ is proximal. We will show that  $x$ is a weakly product recurrent point. For this, let $(Y,S)$ be a minimal system and $y\in Y$. Take $U$ and $V$ be neighborhoods of $x$ and $y$, respectively. We need to show that $N(x,U)\cap N(y, V)\neq\emptyset$.

Since $p_j\ra x$, there is some $p_j\in U$. Take $\epsilon>0$ such that $B(p_j, 2\epsilon)\subset U$. Since $p_j$ is distal, $(p_j,y)$ is a minimal point. Thus $N(p_j, B(p_j, \epsilon))\cap N(y, V)$ is syndetic. On the other hand, $\{n\in \N: \rho(T^nx, T^np_j)<\epsilon\}$ is thick since $(x,p_j)$ is proximal. Thus there is some $n\in N(p_j, B(p_j, \epsilon))\cap N(y, V)$ such that  $ \rho(T^nx, T^np_j)<\epsilon$ and hence
$$T^{n}x\in B(T^{n}p_j,\epsilon)\subset B(p_j,2\epsilon)\subset U.$$ This shows that $N(x,U)\cap N(y, V)\neq\emptyset$. Thus $x$ is weakly product recurrent.
\end{proof}

We do know whether the condition in (1) above is also necessary. We remark here that if the $(X,T)$ in Theorem \ref{wpp with dense distal}  (2) is weakly mixing additionally, then it follows from  \cite{DSY12, Op19} that $(X,T)$ is disjoint from all minimal systems. In this case, any transitive points are weakly product recurrent. However, the converse of Theorem \ref{wpp with dense distal} (2) is not true in general. In \cite{LYY15}, the authors construct a weakly mixing system that is disjoint from all minimal systems but there are no dense distal points. Actually, the system $(X,T)$ they constructed has no dense distal points but the hyperspace $(2^{X},T)$ has dense distal points.

\medskip
Recall that a  {\it $P$-system} is a transitive system with dense periodic orbits.
\begin{cor}Let $(X,T)$ be a transitive $P$-system. Then each $x\in Tran_T$ is a weak product recurrent point.
\end{cor}
\begin{proof}
 Let $(x_i)_{i\in\N}$ be a sequence of periodic points that is dense in $X$. Let $P_i$ be the orbit of $x_i$ with period $n_i$ for each $i\in\N$. Let $k_i=n_1\cdots n_i, i\in \N$.

 Assume that $X_i=\overline{\O}(x, T^{k_i}), i\in \N$. First, we show that $X_i\cap P_j\neq\emptyset$ for every $j\in\{1,\ldots,i\}$. Since $X=X_i\cup TX_i\cup\cdots\cup T^{k_i-1}X_i$, there is some $m\in\{0,1,\ldots,k_i-1\}$ such that $x_j\in T^{m}X_i$. But then $T^{k_i-m}x_j\in T^{k_i}X_i=X_i$. Thus $X_i\cap P_j\neq\emptyset$.

 On the other hand, for  every $j\in\{1,\ldots,i\}$ and every point $y\in X_i\cap P_j$, $(x,y)$ is proximal. Indeed,  for  every point $y\in X_i\cap P_j$, $y$ is a fixed point of $T^{k_i}$ and $y\in \overline{\O}(x, T^{k_i})$. This implies that $(x,y)$ is proximal.

 Clearly, for each $j\in\N$, $X_j\cap P_j\supset X_{j+1}\cap P_j\supset\cdots$. Let $Y=\bigcap_{i=1}^{\infty}X_i$. Then $Y\cap P_j\neq\emptyset$ for each $j\in\N$.
 Thus we can take a sequence $(y_j)$ with $y_j\in Y\cap P_j$ such that $y_j\rightarrow x$. Recall that each $(x,y_j)$ is proximal. Then the conclusion is followed from Theorem \ref{wpp with dense distal}.
\end{proof}

The proof works for $P_i$ is replaced by almost 1-1 extension of adding machine $M_i$ and the fact that $x$ is proximal to some point in $M_i$.
When $(X,T)$ is totaly transitive, the above corollary is not new since it is in $\mathcal{M}^\perp$. When $(X,T)$ is not totaly transitive, it really gives something new, for example $A=\cup \{n_i, n_i+2, \ldots, n_i+2i\}$ is not necessarily a pointwise recurrence sequence.

\medskip
As noted in the introduction, there is little known regarding the properties of a system possessing the weak product property.
Thus, we ask
\begin{ques} Give a characterization of a system possessing the weak product property, particularly for a minimal system with the property.
\end{ques}

\section{$R$-sequences}
In this section we investigate $R$-sequences and prove Theorems \ref{thmb} and \ref{thmc}.

\subsection{Variants of disjointness and the proof of Theorem \ref{thmb}}
Let $(X,T)$ be a transitive system and $(Y,S)$ be a minimal system. Recall that they are weakly disjoint if the product system $(X\times Y, T\times S)$ is transitive. This implies that there is some $x\in Tran_{T}$ and $y\in Y$ such that $\overline{\O}((x,y), T\times S)=X\times Y$. In an equivalent form, there is some $x\in Tran_{T}$ and $y\in Y$ such that for any neighborhood $U$ of $x$,  the orbit of $y$ along $N(x, U)$ is dense in $Y$. A stronger notion than weak disjointness is disjointness. An equivalent characterization of  the disjointness between $(X,T)$ and $(Y,S)$ is that for any $x\in Tran_{T}$ and any $y\in Y$, orbit of $y$ along $N(x, U)$ is dense in $Y$ for any neighborhood of $x$.

Now we introduce a mild notion between weak disjointness and disjointness. For a transitive system $(X,T)$ and a minimal system $(Y,S)$, they are {\it mildly disjoint}, denoted it by $X \rightthreetimes Y$ , if for any $x\in Tran_{T}$ there is some $y\in Y$ such that for any neighborhood $U$ of $x$,  the orbit of $y$ along $N(x, U)$ is dense in $Y$.  Recall that $X\curlywedge Y$ and $X\perp Y$ mean they are weakly disjoint and disjoint, respectively. From the definitions, we immediately have
\[ X\perp Y \Rightarrow X\rightthreetimes Y \Rightarrow X\curlywedge Y.\]
In addition, it is known that disjointness is strictly stronger that weak disjointness. Superficially, mild disjointness is stronger than weak disjointness. To our surprise, they turn out to be equivalent. This is the main result in this subsection and Theorem \ref{thmb} is followed.

\begin{thm}\label{mainsect}
Let $(X,T)$ be a transitive system and $(Y,S)$ be a minimal system. If $X\curlywedge Y$ then $X\rightthreetimes Y$.
\end{thm}

To show Theorem \ref{mainsect}, we start with some lemmas. We suppose that $(X,T)$ is a transitive system and $(Y,S)$ is a minimal system in the sequel.
\begin{lem}\label{0or1} If $(X,T) \curlywedge (Y,S)$, then for any transitive point $x\in X$ either
$D_x=\{y\in Y: (x,y)\in Tran_{T\times S}\}$ is empty or it is a dense  $G_\delta$ set.
\end{lem}

\begin{proof} 
Assume that $D_x\not=\emptyset$. Let $y\in D_x$. Then for any $(z_1,z_2)\in X\times Y$,   choose some point $z_3\in S^{-1}\{z_2\}$ and  there is a sequence $\{n_i\}_{i=1}^\infty\subset \N$ such that $T^{n_i}x\ra z_1, S^{n_i}y\ra z_3\in S^{-1}\{z_2\}$, i.e. $S^{n_i}Sy\ra z_2$. So $Sy\in D_S$. This implies that $D_x$ is dense in $Y$.

Now we show that $D_x$ is residual. For a given $\ep>0$, let $\{V_1,\ldots,V_t\}$ and $\{U_1,\ldots,U_k\}$ be open covers of $X$ and $Y$
respectively with $\diam (V_i), \diam(U_j)<\ep.$ Set
$$D_\ep=\{z\in Y: \forall 1\le i\le t, 1\le j\le k, \exists n\in\N\, \text{such that}\ (T^nx,S^nz)\in V_i\times U_j\}.$$
Then $D_\ep$ is open and dense. It is clear that $D_x=\cap_{n\in\N}D_{1/n}.$
\end{proof}


\begin{lem}\label{F_R}
Let $(X,T)\curlywedge (Y,S)$ and $x\in Tran_T$. Then $D_x\not=\emptyset$ 
if and only if for any open non-empty sets $U$ of $X$ and $V$ of $Y$, $\cup_{n\in N(x,U)}S^{-n}V$ is dense in $Y$.

In fact
$$D_x=\bigcap_{i=1}^\infty\bigcap_{j=1}^\infty (\bigcup_{n\in N(x,U_i)}S^{-n}V_j),$$
where $\{U_i\}$ and $\{V_i\}$ are bases of the topologies of $X$ and $Y$ respectively.

\end{lem}
\begin{proof} 
First we assume that there is $y\in Y$ such that $\overline{\O((x,y),T\times S)}=X\times Y$. We claim that $\{S^ny: n\in N(x,U)\}$ is dense in $Y$ for any nonempty open subset $U$ of $X$.

For any $z\in Y$ and any neighborhood $U_z$ of $z$, there is $n\in \N$ such that $(T^nx,S^ny)\in U\times U_{z}$. It is clear that $n\in N(x,U)$. It implies that there is $n\in N(x,U)$ such that $S^ny\in U_z$. This means that $\{S^ny: n\in N(x,U)\}$ is dense in $Y$. By Lemma \ref{s=dence} we get that $\cup_{n\in N(x,U)}S^{-n}V$ is dense in $Y$ for any open non-empty set $V$.

\medskip
Now we show the other direction. We need to show there is $y\in Y$ with $(x,y)\in Tran_{T\times S}$.

By the assumption for each non-empty open set $V$ of $Y$, $\cup_{n\in N(x,U)}S^{-n}V$ is dense in $Y$ for open non-empty set $U$ of $X$. Fix a transitive point $z\in U$ and neighborhood $U_i\subset U$ of $z$ with $\diam(U_i)\ra 0$.
Set $$\Omega=\bigcap_{i=1}^\infty\bigcap_{j=1}^\infty (\bigcup_{n\in N(x,U_i)}S^{-n}V_j),$$ where $\{V_j\}$ is a base of the topology of $Y$. It is clear that $\Omega$ is a dense $G_\delta$ set of $Y$.

Let $y\in \Omega$. Then for each $i,j\in\N$ there is $n\in N(x,U_i)$ such that $T^nx\in U_i$ and $S^ny\in V_j$. This implies that $\overline{\O((x,y),T\times S)}\supset \{z\}\times Y$, and hence $\overline{\O((x,y),T\times S)}=X\times Y$ since $z$ is a transitive point.
\end{proof}


The following lemma is crucial for our proof.
\begin{lem}\label{dense-orbit-set}
Let $(X,T)$ be a minimal system and $A\subset \N$. If there is some $x\in X$ such that ${orb}_{A}(x)$ is dense in $X$, then for any nonempty open sets $U,V\subset X$, there is a finite subset $B\subset A$ such that for any $k\in\mathbb{N}$, $\left(\bigcup_{n\in B}T^{-n}U\right)\cap T^{-k}V\neq\emptyset.$
\end{lem}
\begin{proof}
Fix nonempty open sets $U,V$ in $X$. Since $(X,T)$ is minimal, there is some $r\in\mathbb{N}$ such that
\[X=V\cup T^{-1}V\cup\cdots\cup T^{-r}V.\]

\noindent{\bf Claim}. There is some $y\in X$ and a finite subset $B\subset A$ such that
 \[\{y, Ty, \ldots, T^{r}y\}\subset \bigcup_{n\in B}T^{-n}U.\]

First, we show that the Claim implies the proposition. For any $k\in\mathbb{Z}$, there is some $i\in\{0,1,\ldots,r\}$ such that $T^{k}y\in T^{-i}V$. Thus $T^{i}y\in T^{-k}V$. This implies that $\{y, Ty, \ldots, T^{r}y\}\cap T^{-k}V\neq\emptyset$ for any $k\in\N$. In particular, $ (\bigcup_{n\in B}T^{-n}U)\cap T^{-k}V\neq\emptyset$ for any $k\in\N$.

Now we are going to show the Claim. By the assumption and Lemma \ref{s=dence}, $\bigcup_{n\in A}T^{-n}U$ is dense in $X$. Thus there is some $n_0\in A$ such that $U\cap T^{-n_0}U\neq\emptyset.$ By the same reason, there is some $n_1\in A$ such that
\[ (U\cap T^{-n_0}U) \cap (T^{-n_1}T^{-1}U)\neq\emptyset.\]
Continuing this process, we can find $n_0,n_1,\ldots,n_r\in A$ such that
\[ W:=(U\cap T^{-n_0}U) \cap (T^{-n_1}T^{-1}U)\cap \cdots\cap (T^{-n_r}T^{-r}U)\neq\emptyset.\]
Take $y\in W$. Then $y\in U, Ty\in T^{-n_1}U, T^{2}y\in T^{-n_2}U,\ldots, T^{r}y\in T^{-n_r}U$. Set $B=\{n_0,n_1,\ldots,n_r\}$. Then we have
 \[\{y, Ty, \ldots, T^{r}y\}\subset \bigcup_{n\in B}T^{-n}U.\]

\end{proof}

\begin{proof}[Proof of Theorem \ref{mainsect}]
Let $x\in Tran_T$ and $y'\in Y$ such that $(x,y')\in Tran_{T\times S}$. By the proof of the first part of  Lemma \ref{F_R},
$\{S^ny':\ n\in N(x,U)\}$ is dense in $Y$ for any open non-empty set $U$ of $X$.


Let $x'\in Tran_T$.  Further, fix nonempty open sets $U\subset X$ and $V,W\subset Y$. By Lemma \ref{F_R}, it suffices to show that there is some $n\in N(x', U)$ such that $(S^{-n}V)\cap W\neq \emptyset$.

By Lemma \ref{dense-orbit-set}, there is a finite set $F\subset N(x, U)$ such that for any $k\in \mathbb{Z}$,
\begin{equation}\label{eq1}
 \left(\bigcup_{m\in F}S^{-m}V \right)\cap S^{-k} W\neq\emptyset.
 \end{equation}
Since $x'$ is a transitive point, it follows from the continuity that there is some $\ell\in \mathbb{Z}$ such that $\{ T^{m+\ell}x': m\in F\}\subset U$, since $T^mx\in U$ for any $m\in F$ and there is a sequence $\{l_i\}\subset \N$ such that $T^{l_i}x'\ra x$. In other words, $F+\ell\subset N(x', U)$. According to (\ref{eq1}), there is some $m\in F$ such that
\[ (S^{-(m+\ell)} V)\cap W\neq \emptyset.\]
This implies that there is some $n\in N(x', U)$ such that $(S^{-n}V)\cap W\neq \emptyset$, since $n=:m+\ell\in N(x', U)$.
\end{proof}

\subsection{Scattering and the proof of Theorem \ref{thmc}}

 Recall that a system is {\it scattering} if it is weakly disjoint from every minimal system and is {\it weakly scattering} if it is weakly disjoint from every minimal equicontinuous system. In the second case, it is equivalent to be disjoint from every minimal equicontinuous system.

Akin and Kolyada \cite{AG03} showed that if $(X,T)$ is weakly mixing then for each $x\in X$, $P[x]$ is a dense $G_\delta$-set. The following corollary-(4) improves the result for a minimal weakly mixing system.

\begin{cor}\label{manycor} We have
\begin{enumerate}
\item If  $(X,T)$ is scattering, then every $x\in Tran_T$ is a scattering point.

\item Assume that $(X,T)$ is transitive, $(Y,S)$ is minimal and $X\curlywedge Y$. Then for each $x\in Tran_T$, $D_x=D_x(Y)$ is a dense $G_\delta$ set of $Y$, where $D_x=\{y\in Y: (x,y)\in Tran_{T\times S}\}$.
\item Let $(X,T)$ and $(Y,S)$ be minimal and $X\curlywedge Y$. Then for each $x\in X$, $D_x$ is a dense $G_\delta$ set of $Y$.
\item Assume that $(X,T)$ is scattering and $(Y,S)$ is  an $M$-system. Then for each $x\in Tran_T$, $D_x$ is a dense $G_\delta$ set of $Y$.

 Particularly, if $(X,T)$ is a weakly mixing, minimal system (resp. $M$-system), then for each $x\in X$ (resp. $x\in Tran_T$), $D_x$ is a dense $G_\delta$ set of $Y$.
\end{enumerate}
\end{cor}
\begin{proof} (1), (2) and (3)  follow from Theorem \ref{mainsect} and Lemma \ref{0or1}.

\medskip
Now we show (4). Let $\{M_i\}\subset Y$ be minimal subsets whose union is dense in $Y$.
Let $U,V$ be open non-empty sets of $X$ and $Y$ respectively. By Theorem \ref{mainsect}, $D_x(M_i)\not=\emptyset$. By Lemma \ref{F_R},
 $\cup_{n\in N(x,U)}S^{-n}(V\cap M_i)$ is dense in $M_i$ if $M_i\cap V\not=\emptyset$, since $X\curlywedge M_i$.

Let $I=\{i\in \N: M_i\cap V\not=\emptyset\}$. We claim that $\cup_{i\in I}M_i$ is dense in $Y$. Assume the contrary that $\overline{\cup_{i\in I}M_i}\not=Y.$ Then there is a non-empty open set $V'\subset Y\setminus \overline{\cup_{i\in I}M_i}$. The transitivity of $(X,T)$ implies that
there is $M_j$ such that $M_j\cap V', M_j\cap V\not=\emptyset$, a contradiction. This proves the claim.

Thus, we have that  $\cup_{n\in N(x,U)}S^{-n}V$ is dense in $Y$. Again using Lemma \ref{F_R}, it follows that $D_x(Y)\not=\emptyset$. Lemma \ref{0or1} implies that $D_x$ is a dense $G_\delta$-set.

Since weakly mixing implies scattering, the argument applies to the case when $(X,T)$ is a  weakly mixing $M$-system system.
\end{proof}

\begin{ques}
Is there a weakly mixing system $(X,T)$ satisfying that there is $x\in Tran_T$ with $(x,y)\not\in Tran_{T\times T}$
for any $y\in X$? 
\end{ques}

We remark that such a system is not an $M$-system by Corollary \ref{manycor}(4). So, we may consider to find counterexamples in the class of
proximal weakly mixing systems.

\medskip

It is known that if $(X,T)$ is an $E$-system, then weakly scattering implies weak mixing (see Lemma \ref{E-w-s}). We have

\begin{thm}\label{scatterin-char-w}
 Let $(X,T)$ be an $M$-system. Then the following statements are equivalent.
\begin{enumerate}
\item $(X,T)$ is weakly mixing.
\item $(X,T)$ is scattering.
\item $(X,T)$ is weakly disjoint from all minimal subsystems of $(X,T)$.
\item There are minimal sets $M_i,i\in\N$ whose union is dense in $X$ and $(X,T)$ is weakly disjoint from each $M_i$.
\item There are minimal sets $M_i,i\in\N$ with $\overline{\cup_{i\in\N}M_i}=X$ such that for each $x\in Tran_T$ and for each $i\in\N$, $x$
 is proximal to each point in a dense $G_\delta$ set $M_{i,x}$ of $M_i$.
 \item There are minimal sets $M_i,i\in\N$ with $\overline{\cup_{i\in\N}M_i}=X$ such that for some $x\in Tran_T$ and for each $i\in\N$, $x$
 is proximal to each point in a dense invariant set $M_{i,x}$ of $M_i$.

\end{enumerate}
\end{thm}
\begin{proof} It is clear that $(i)$ implies $(i+1)$ for $i=1,2,3$.


Now we show (4) implies (5). Fix $i\in \N$ and $x\in Tran_T$. Since $X\curlywedge M_i$, by Theorem \ref{mainsect} there is $z\in M_i$ such that
$(x,z)$ is a transitive point of $X\times M_i$. Then by Lemma \ref{0or1} the set $M_{i,x}=\{z\in M_i: (x,z)\in Tran_{X\times M_i}\}$ is a dense
$G_\delta$ set. It is clear that $(x,z)$ is proximal for each $z\in M_{i,x}$, since $\overline{\O}((x,z),T\times T)\supset \Delta_{M_i}$.

It is clear (5) implies (6).

\medskip

(6) implies (1). Let $i\in \N$. Pick $y\in M_{i,x}$ and consider $y'=T^ny$ for $n\in \N$.
Then $(x,y)$ is proximal and $(x,y')$ is proximal, since $M_{i,x}$ is invariant. Note that $(y,y')$ is minimal. So for any $\ep>0$,
$\{n\in \N: \rho(T^nx,T^ny'\}\in \F_t$ and $\{n\in\N: \rho(T^ny',y')<\ep,\rho(T^ny,y)<\ep\}\in \F_s$.
This shows that $(y',y)\in \overline{\O}((x,y), T\times T)$ and so $M_i\times \{y\}\subset \overline{\O}((x,y),T\times T)$ which implies that $M_i\times  M_i\subset \overline{\O}((x,y),T\times T)$ by the minimality of $M_i$.

Now let $U_1,\ldots,U_4$ be open non-empty sets of $X$. As $x\in Tran_T$ there are $n_1,\ldots,n_4\in \N$ such that $T^{n_i}x\in U_i$, $1\le i\le 4$. Assume that $x_i\in M_i$ with $\lim_{i\ra \infty} x_i=x$. The continuity of $T$ implies that there is $i_0$ such that  $T^{n_i}x_{i_0}\in U_i$,
$1\le i\le 4$. Let
$M$ be the orbit closure of $x_{i_0}$. Then $M\cap U_j\not=\emptyset$, $1\le j\le 4$. So, there are $y\in M$, $n\not=m\in\N$ such that
$$(T^nx,T^ny)\in U_1\times U_2 \  \text{and}\ (T^mx,T^my)\in U_3\times U_4.$$ Thus, $x\in T^{-n}U_1\cap T^{-m}U_3$ and $y\in T^{-n}U_2\cap T^{-m}U_4$ which implies that $X$ is weakly mixing.

\end{proof}

\subsection{Recurrence sequences arising from disjointness}
Recall that $A\subset \mathbb{N}$ is an $R$-{\it sequence} if for any any minimal system $(X,T)$, there is some $x\in X$ such that the orbit of $x$ along $A$ is dense in $X$, that is $orb_{A}(x):=\{T^{n}x:\ n\in A\}$ is dense. This is closely related to the notion of dense orbit set introduced in \cite{GTWZ21}. A subset set $A\subset \N$ is a {\it dense orbit set} if for any minimal system $(X,T)$ and any $x\in X$, the orbit of $x$ along $A$ is dense in $X$.

The collections of dense orbit sets is  denoted by $\mathcal{F}_{DOS}$. 

The first statement in the following Theorem is given in \cite{HY05}.
\begin{thm}\label{equi-s-point}
Let $(X,T)$ be a transitive system. Then
\begin{enumerate}
\item[(1)] $X\perp \mathcal{M}$ if and only if for any $x\in Tran_{T}$ and any neighborhood $U$ of $x$, $N(x, U)$ is a dense orbit set;
\item[(2)] $X\curlywedge \mathcal{M}$ if and only if for any $x\in Tran_{T}$ and any neighborhood $U$ of $x$, $N(x, U)$ is an $R$-sequence.
\end{enumerate}
\end{thm}
\begin{proof}We show the second assertion.

Suppose that $X\curlywedge\mathcal{M}$ and take $x\in Tran_{T}$.  Let $(Y,S)$ be a minimal system. By Theorem \ref{mainsect}, we have $X\rightthreetimes Y$. Thus there is some $y\in Y$ such that $(x,y)\in Tran_{T\times S}$. In particular, $\{S^{n}y:\ n\in N(x,U)\}$ is dense in $Y$ for any neighborhood $U$ of $x$. Thus $N(x, U)$ is an $R$-sequence.

Now suppose that $x\in Tran_{T}$ and $N(x, U)$ is an $R$-sequence,  for any neighborhood $U$ of $x$. Let  $(Y,S)$ be a minimal system. We need to show that $X\curlywedge Y$. To this end, we take a countable neighborhood basis $(U_{i})_{i\in\mathbb{N}}$ of $x$. Then for each $i\in\mathbb{N}$, there is some $y_i\in Y$ such that $\{S^{n}y:\ n\in N(x,U_i)\}$ is dense in $Y$. Further, it follows from Lemma \ref{s=dence} that
\[ Y_{i}:=\left\{y\in Y: \overline{\{S^{n}y:\ n\in N(x,U_i)\}}=Y\right\}\]
is a dense $G_{\delta}$ set. Then $\Omega:=\bigcap_{i=1}^{\infty}Y_i$ is also a dense $G_{\delta}$ set. Now we can claim that $\{S^{n}y:\ n\in N(x,U)\}$ is dense in $Y$ for any $y\in \Omega$ and any neighborhood $U$ of $x$. Indeed, there is some $U_i\subset U$ and hence $\{S^{n}y:\ n\in N(x,U)\}\supset \{S^{n}y:\ n\in N(x,U_i)\}$. Thus  $\{S^{n}y:\ n\in N(x,U)\}$ is dense in $Y$. Therefore, $\O((x,y), T\times S)$ is dense in $X\times Y$ since $x\in Tran_{T}$. This shows that $X\curlywedge Y$ and hence $X\curlywedge\mathcal{M}$.
\end{proof}

Since disjointness is strictly stronger than weak disjointness,  it is reasonable to conclude that $\mathcal{F}_{R}$ is strictly larger than $\mathcal{F}_{DOS}$.   It is known that a dense orbit set must be piecewise syndetic (\cite{HY05}). In other words, $\mathcal{F}_{DOS}\subset\mathcal{F}_{ps}$. Next, we give an example of an $R$-sequence that is not piecewise syndetic.

\begin{ex}\label{ex-example1} There is $A\in \F_R$  such that $A\not\in \F_{ps}.$
\end{ex}
\begin{proof} Let $\{t_i\}$ be a subsequence of $\N$ such that for each $n\in \N$, $\{i\in\N:t_i=n\}$ is infinite, saying for example
$(t_i)=(1,2,1,2,3,1,2,3,4,\ldots)$.

Let $A_i$ be a sequence of finite IP-sets with $|A_i|\ra \infty$ and
$$\max (A_i+t_i)<\min (A_{i+1}+t_{i+1})$$ for each $i\in \N$.
Then $A=\cup_{i=1}^\infty (A_i+t_i)$ is the set we want.

By Theorem \ref{rdop-char}, it remains to show that  for a minimal system $(X,T)$ and an open non-empty subset $W$ and $n\in \N$, $A\cap (n+N(W,W))\not=\emptyset$.

Assume that $\mu$ is an ergodic measure on $X$. Take $i\in \N$ such that $A_i$ is generated by $\{n_1<\cdots<n_k\}$ wit $k>[\frac{1}{\mu(W)}]+1$ and $t_i=n$. Consider
$$\{T^tW: t=n_1, n_1+n_2,\ldots, n_1+\ldots+n_k\}.$$ It is clear that $\mu(T^lW)=\mu(W)>0$ for each $l\in \N$.
So there are $1\le k_1<k_2\le k$ such that $\mu(T^{n_1+\cdots+n_{k_1}}W\cap T^{n_1+\cdots+n_{k_2}}W)>0$, i.e.
$m=:n_{k_1+1}+\cdots+n_{k_2}\in N(W,W)$. Since $m+n\in A$, it implies that
$A\cap (n+N(W,W))\not=\emptyset$.
\end{proof}

According to Theorem \ref{F_R} we have many examples of $\F_R$.

\begin{ex} Let $(X,T)$ be a scattering system. Then for any any transitive point $x$ and non-empty open subset $U$ of $X$, $N(x,U)\in \F_R$.
\end{ex}

Since weakly mixing systems are scattering, we have

\begin{ex} Let $(X,T)$ be a weakly mixing system. Then for each $x\in X$ and each open non-empty set $U$, $N(x,U)\in \F_R$.
\end{ex}

In addition, we also have the following example to give an $R$-sequence that is not a dense orbit set.
\begin{ex}\label{F-R-example} Let $(X,T)$ be a scattering almost equicontinuous system.  Then each transitive point $x$ of $(X,T)$ is a scattering point and there is a neighborhood $U$ of $x$ such that $N(x,U)\not\in \F_{DOS}$. Particularly, the system in Lemma \ref{AG-exa} is a such system.
\end{ex}

\subsection{Katznelson and Richter's Questions}
 A set $A\subset \N$ is called a {\it Bohr recurrence set} if for every rotation on a torus of finite dimension  $(\T^{n}, T)$ there is some $x\in \T^{n}$ and subsequence $(n_i)$ in $A$ such that $T^{n_i}x\to x$ as $i$ tends to the infinity. A longstanding open question in topological dynamics is Katznelson's question, saying whether topological recurrence sets and Bohr recurrence sets coincide.  In \cite{HSY16}, Huang, Shao and Ye raised a higher order version of Katznelson's question, saying that whether multiply recurrence sets are equivalent to Nil-Bohr recurrence sets.  Recently, Alweiss gives a negative answer to this higher order version of Katznelson's question in the preprint \cite{Al25}. This suggests that Katznelson's question may also have a negative answer.

  It was proved in \cite[Theorem 4.10]{HY02} that $(X,T)$ is scattering if and only if for any nonempty open sets $U,V$ in $X$,  $N(U,V)$ is a topological recurrence sequence, and in \cite[Theorem 2.11]{HY12} weakly scattering if and only if $N(U,V)$ is a Bohr$_0$ set. Another open question  is whether there is a weak scattering system that is not scattering. In \cite{HY12}, Huang and Ye show that the existence of such a system will answer Katznelson's  question in the negative. By \cite[Theorem 2.14]{HY12} a weakly scattering system is disjoint from all minimal equicontinuous systems and if $(X,T)$ is transitive and $X\perp\mathcal{M}$, then it is weakly mixing \cite[Theorem 2.6]{HY05}. Note that every weakly mixing system is scattering. In \cite{HY02}, Huang and Ye constructed a scattering system that is not weakly mixing. In particular, it is not disjoint from all minimal systems.

Katznelson's question is asked whether topological recurrence can be captured by the recurrence in equicontinuous systems, while the equivalence between scattering and weak scattering property is asked whether weak disjointness from all minimal systems can be captured by the weak disjointness from all minimal equicontinuous systems. In the private communication, Richter asked a question in a similar favor \cite[Question 17]{HSQY26}.

 Recall that $A\subset \N$ is an $R_e$-{\it sequence} if  for each minimal equicontinuous system $(X,T)$ there is $x\in X$ such that $\{T^nx:\ n\in A\}$ is dense in $X$. The collection of $R_e$-sequences is denoted by $\F_{R_e}$. It is clear that $\F_R\subset \F_{R_e}$. Richter asked the following question.
\begin{ques}\label{FR=FER?}
Is it true that $\F_R=\F_{R_e}$?
\end{ques}

By the same proof in Theorem \ref{equi-s-point} we have

\begin{thm}\label{F-ER} Let $(X,T)$ be transitive and $x$ be a transitive point. Then for each open non-empty set $U$ of $x$, $N(x,U)\in \F_{R_e}$
if and only if for each minimal equicontinuous system $(Y,S)$ there is $y\in Y$ such that $\overline{\O}((x,y), T\times S)=X\times Y$.
\end{thm}
We remark that the existence of a point $y\in Y$ with $\overline{\O}((x,y), T\times S)=X\times Y$ above implies  that  $\overline{\O}((x,y'), T\times S)=X\times Y$ for every point $y'\in Y$ due to the equicontinuity.

 Next, we show that a negative answer to Richter's question also leads to a negative answer to Kazhnelson's question.

 \begin{thm}\label{Ritcher+Katznelson}
 \begin{itemize}
\item[(1)] If there is a weakly scattering system that is not scattering, then $\mathcal{F}_{R}\neq \mathcal{F}_{R_e}$.
\item[(2)] If $\mathcal{F}_{R}\neq \mathcal{F}_{R_e}$ then there is a Bohr recurrence set that is not a topological recurrence set.
 \end{itemize}
 \end{thm}
\begin{proof}
(1) Suppose that $(X,T)$ is weakly scattering but not scattering. Then there is a minimal system $(Y,S)$ such that $(X\times Y, T\times S)$ is not transitive. In particular, there is $x\in Tran_{T}$, a neighborhood $U$ of $x$ and a nonempty open set $V$ of $Y$ such that $\bigcup_{n\in N(x,U)}T^{-n}V$ is not dense in $Y$. This implies that $N(x,U)\notin\mathcal{F}_{R}$. However, $N(x,U)\in\mathcal{F}_{R_e}$, since $X$ is disjoint from every minimal equicontinuous system.

(2) Suppose that $\mathcal{F}_{R}\neq \mathcal{F}_{R_e}$ and we take $A\in\mathcal{F}_{R_e}\setminus\mathcal{F}_{R}$. Then there is a minimal system $(X,T)$ such that for any $x\in X$, $\{T^{a}x:\ a\in A\}$ is not dense in $X$. We claim that there is some $k\in\N$ and a nonempty open set $U\subset X$ such that $(\bigcup_{a\in A}T^{-(a+k)}U)\cap U=\emptyset$.  This implies that $A+k$ is not a topological recurrence set. By assumption, there is a nonempty open set $V$ of $X$ such that $\bigcup_{a\in A}T^{-a}V$ is not dense in $X$. Take $x\in V$. By the minimality of $(X,T )$, there is some $k\in \N$ such that $T^kx\not\in \overline{\bigcup_{a\in A}T^{-a}V}$. Further, there is some  neighborhood $W$ of $x$  such that $(\bigcup_{a\in A}T^{-a}V)\cap T^{k}W=\emptyset$. Let $U=V\cap W$ and then $(\bigcup_{a\in A}T^{-a}U)\cap T^{k}U=\emptyset$ and hence $(\bigcup_{a\in A}T^{-(a+k)}U)\cap U=\emptyset$. Since $A\in \mathcal{F}_{R_e}$, it is easy to see that $A+k\in \mathcal{F}_{R_e}$. In particular, $A+k$ is a Bohr recurrence set. This completes the proof.
\end{proof}

Next, we discuss some positive aspects towards Katznelson's question. Although we do not know whether a weakly scattering system is weakly disjoint from any minimal system. it is known that a weakly scattering system is indeed disjoint from every minimal distal system (see \cite{HY12}). However, it is still open that whether a Bohr recurrence set is also a recurrence set for all distal systems. In \cite{GKR22}, Glasscock et show that it is the case for some special skew product systems. So it is natural to expect that $\mathcal{F}_{R_e}=\mathcal{F}_{R_d}$, where  $\mathcal{F}_{R_d}$ is the collection of subset $A\subset \N$ such that for each minimal distal system contains a point the orbit of which under $A$ is dense.  At least, we hope the following question has a positive answer.

\begin{ques}
If $A\in\mathcal{F}_{R_e}$ then any minimal nilsystem  contains a point the orbit of which under $A$ is dense.
\end{ques}

\section{$R$-sequences for totally minimal systems}

In ergodic theory it was shown for certain sequences $A=\{a_1,a_2,\ldots\}$ (for example, polynomial sequences and prime numbers), the ergodic average  $\frac{1}{N}\sum_{i=1}^{N}T^{a_i}f$ is captured by the ergodic average in the rational factor. In particular, in such cases we have for a totally ergodic (TE for short) measure preserving system $(X,\mathcal{B},\mu,T)$ and $f\in L^2(\mu)$,
\begin{equation}\label{strong-p}
\frac{1}{n}\sum_{i=1}^n f(T^{a_i}x)\ra \int_X fd\mu.
\end{equation}

So, it is natural to consider the topological counterparts. That is, we are interested in sequences $A=\{a_1,a_2,\ldots\}$ such that
for any totally minimal (TM for short) system $(X,T)$ there is $x\in X$ such that the orbit of $x$ along $A$ is dense in $X$.

\medskip

Since a totally minimal system does not necessarily carry a totally ergodic measure, it is not directly that the sequence satisfying (\ref{strong-p}) for TE is  an $R$ sequence for TM. Philosophically, one thinks that it is true. For example, one conjectures that the set of primes is an $R$-sequence for TM.

\subsection{$\F_{RTM}$}
A set $A\subset \N$ is an $R$-sequence for TM if for any totally minimal system $(X, T)$, there is $x\in X$ such that the set $\{T^nx: n \in A\}$ is dense in $X$. The collection of all such sequences is denoted by $\F_{RTM}.$
It is clear that $\F_R\subset \F_{RTM}$.

\medskip
It is known that $f$ is an integral polynomial with $f(0)=0$ and $f(\N)\subset \N$, then $A=\{f(n):n\in\N\}\in \F_{RTM}$, see \cite{GHSWY20, QJ23}. We now show that $\{n^2: n\in\N\}\in \F_{RTM}\setminus \F_R.$

Let $X=\{x_1,x_2,x_3\}$ be a periodic orbit of period 3 and $T:X\ra X$ such that $T(x_i)=x_{i+1(\text{mod}\ 3)}$ for any $1\le i\le 3$. Since  $n^2\equiv 1\mod 3$ for any $n\in\N$,  we conclude that $\{T^{n^2}x_i:n\in\N\}$ is not dense in $X$ for any $i\in\{1,2,3\}$.

 Next, we show that $\{f(n):n\in\mathbb{N}\}\not\in\mathcal{F}_{R}$ whenever $\deg f\geq 2$. To this end, we need some results on permutational polynomials (\cite[Chapter 8]{MP13}).

  Let $p$ be a prime and $q=p^{r}$. Let $\mathbb{F}_{q}$ be a finite field with $q$ elements.

\begin{defn}
A polynomial $f\in \mathbb{F}_{q}[x]$ is called a {\it permutational polynomial} of $\mathbb{F}_{q}$ if $f: c\mapsto f(c)$ is a permutation of $\mathbb{F}_{q}$.
\end{defn}

\begin{thm}[Hermite's Criterion]
A polynomial $f\in \mathbb{F}_{q}[x]$ is a permutational polynomial of $\mathbb{F}_{q}$ if and only if the following conditions hold:
\begin{enumerate}
\item the reduction of $f(x)^{q-1}\mod (x^{q}-x)$ is monic of degree $q-1$;
\item for each integer $k$ with $1\leq k\leq q-2$ and $k\not\equiv {\rm mod } p$, the reduction  of $f(x)^{k}\mod (x^{q}-x)$ has degree $\leq q-2$.
\end{enumerate}
\end{thm}

\begin{cor}\label{xuhui}
Let $f(x)\in \mathbb{Z}[x]$. If $\deg f(x)\geq 2$, then there is some prime $p$ such that
\[ \{f(n) ~{\rm mod }~p: n\in\mathbb{Z}\}\neq \{0,1,\ldots, p-1\}.\]
\end{cor}
\begin{proof}
Let $f(x)=c_{a}x^{a}+\cdots+c_1x+c_0$ with $a\geq 2$ and $c_a\neq 0$. Then there are infinitely many prime numbers of the form $an+1$. Take such a prime $p=an+1$ such that $p> |c_a|$.  Assume that $\{f(n) ~{\rm mod }~p: n\in\mathbb{Z}\}= \{0,1,\ldots, p-1\}$.

Let $\overline{f(x)}$ be the homomorphism image of $f(x)$ in $\mathbb{F}_{p}[x]$. Then $\overline{f(x)}$ is a permutational polynomial of $\mathbb{F}_{p}$. Since $p> |c_a|$,
$\deg \overline{f(x)}=\deg f(x)=a$. Further,
\[\deg \overline{f(x)}^{n}=an=p-1.\]
Since $1\leq n\leq p-2$, it follows from the Hermite's criterion that $\overline{f(x)}$ is not a permutational polynomial of $\mathbb{F}_{p}$. Hence
\[ \{f(n) ~{\rm mod }~p: n\in\mathbb{Z}\}\neq \{0,1,\ldots, p-1\}.\]
\end{proof}

Now let $(X,T)$ be a periodic orbit of period $p$, where $p$ is a prime number defined in Corollary \ref{xuhui}.
Then for any $x\in X$, $\{T^{f(n)}x:\ n\in\N\}$ is not dense in $X$. Thus $\{f(n):\ n\in\N\}\notin \mathcal{F}_{R}$.

\medskip

Similarly to Theorem \ref{rdop-char} we have
\begin{thm}\label{RTM-char}
For $A\subset \mathbb{N}$, the following statements are equivalent
\begin{enumerate}
\item $A\in \F_{RTM}$.
\item  for any totally minimal system $(X,T)$, any open non-empty sets $U$ and $V$ of $X$, $A\cap N(U,V)\not=\emptyset$, i.e. there is $n\in A$ such that $U\cap T^{-n}V\not=\emptyset$.
\item  for any totally minimal system $(X,T)$, any open non-empty set $U$ of $X$ and $n\in\Z$, $A\cap (n+N(U,U))\not=\emptyset$.
\item for any totally minimal system $(X,T)$, any open non-empty set $U$ of $X$, $\cup_{a\in A}T^{-a}U$ is dense in $X$.
\item for any totally minimal system $(X,T)$ there is a dense $G_\delta$ set $\Omega$ of $X$ such that for each $x\in \Omega$,
$\{T^ax:a\in A\}$ is dense in $X$.
\end{enumerate}
\end{thm}



\subsection{Weak disjointness from all totally minimal systems} We may carry all in the discussions as in the previous  section. That is,
we can consider systems weakly disjoint (resp.  disjoint) for TM.


In \cite{HY05} it was proved that if $(X,T)$ is transitive and $x$ is a transitive point, then $(X,T)$ is disjoint from all minimal systems if and only if for each neighborhood $U$ of $x$, $N(x,U)\cap A\not=\emptyset$, where $A$ is any $m$-set. Here, $A$ is an $m$-set (termed a dynamical syndetic set in \cite{GL24}) if there is a minimal system $(Y,S)$ with $y\in Y$ and a nonempty open subset $V$ of $Y$ such that $A\supset N(y,V)$ (see \cite{HY05}).

So we have the following theorem by the same argument.
\begin{thm}\label{rec char of dis TM}
 Let $(X,T)$ be transitive and $x$ be a transitive point. Then
\begin{enumerate}
\item[(1)] $(X,T)$ is disjoint from all TM systems
if and only if for each open non-empty subset $U$ of $X$, $N(x,U)\cap A\not=\emptyset$, where $A$ is any m-set for TM.
\item[(2)] $(X,T)$ is disjoint from all TM systems
if and only if for each open non-empty subset $U$ of $X$, $N(x,U)\in \mathcal{F}_{RTM}$.
\end{enumerate}
\end{thm}
Where an $m$-set for TM means that there is a TM system $(Y,S)$, $y\in Y$ and a non-empty open set $U\subset Y$ such that
$A=N(y,U)$.

\medskip

The following lemma is a collection of some known results.

\begin{lem} \label{weaklydsijointmm} Let $X,Y$ be minimal. Then
\begin{enumerate}
\item $X\curlywedge Y$ if and only if $X_{eq}\perp Y_{eq}$ (\cite[Corollary 6]{Peleg}).
\item If $(X,T)$ is distal then $X\perp Y$ if and only if $X\perp Y_{eq}$ (\cite[Chapter 11, Theorem 7,Theorem 9]{Aus}).
\item If a minimal $X'$ is a proximal extension of $X$ and $X\perp Y$, then $X'\perp Y$ (see \cite[Chapter 11, Corollary 13]{Aus}).

\item If a minimal $X'$ is a distal extension of $X$ and $X'\times Y$ is transitive then $X\perp Y$ if and only if $X'\perp Y$
(see \cite[Proposition 7.5]{DSY12})

\item Let $(X,T)$ an $E$-system and $(Y,S)$ be minimal. If $X\perp Y_{eq}$ then $X\curlywedge Y$ (\cite[Proposition A.1]{BHM})
\end{enumerate}
\end{lem}

Unlike the case of $\mathcal{M}^\perp$, it is possible that there is a non-trivial minimal system which is in $\mathcal{M}^\perp_{TM}$.
To see this, first note that a periodic orbit has the property and thus the inverse limit of periodic orbits. This implies that all adding machines
has the property. Generally we have following result.

First, we recall the notion of rational Kronecker factor of a minimal system. Let $(X,T)$ be a minimal system. A continuous function $f$ on $X$ with $|f|=1$ is an {\it eigenfunction} of $T$ if there is a nonzero $\lambda\in\mathbb{C}$ such that $f(Tx)=\lambda f(x)$ for every $x\in X$. In this case, $\lambda$ is called an {\it eigenvalue} of $T$. Now let $\Lambda$ be the collection of eigenvalues of $T$. For each $\lambda\in \Lambda$, let $f_{\lambda}$ be the corresponding eigenfunction. Let $\widehat{\Lambda}$ be the Pontrjagin dual of $\Lambda$. Further, let $\alpha$ be the inclusion map from $\Lambda$ to the unit circle $\mathbb{S}^{1}$. Then $\alpha$ is a character of $\Lambda$ and we identify it with an element of $\widehat{\Lambda}$. Note that for each $x\in X$, $\Lambda\rightarrow \mathbb{S}^{1}, \lambda\mapsto f_{\lambda}(x)$ is a character of $\Lambda$. One can verify that
\[\pi: X\rightarrow \widehat{\Lambda},\ x\mapsto f_{\lambda}(x)\]
is a factor map between $(X, T)$ and $(\widehat{\Lambda}, R_{\alpha})$, where $R_{\alpha}: \widehat{\Lambda}\rightarrow\widehat{\Lambda}, z\mapsto z+\alpha$. Moreover, $(\widehat{\Lambda}, R_{\alpha})$ is isomorphic to the maximal equicontinuous factor of $(X,T)$ (See \cite[Chapter 3]{HK18} for more details).

Let $\Lambda_{rat}=\{\lambda\in \Lambda: \ \exists n\in\N \text{ such that } \lambda^{n}=1\}$ be the set of rational eigenvalues. Then $\widehat{\Lambda_{rat}}$ is a quotient of $\widehat{\Lambda}$ and the induced factor $(\widehat{\Lambda_{rat}}, R_{\alpha})$ is  the {\it maximal rational Kronecker factor} of $(X,T)$. We also use $X_{rat}$ to denote the maximal rational Kronecker factor of $(X,T)$.

Further, for each $\lambda\in\Lambda$, let $\Lambda_{\lambda}$ be the subgroup generated by $\lambda$. Then $(\widehat{\Lambda_{\lambda}}, R_{\alpha})$ is also a factor of $(X,T)$.

\begin{thm}\label{weakdisjointtmminimal}
Let $(X,T)$ be minimal. Then
\begin{enumerate}
\item $X\curlywedge \mathcal{M}_{TM}$ if and only if $X_{eq}=X_{rat}$, 
    where $X_{rat}$ is the maximal rational Kronecker factor of $(X,T)$.
\item If $(X,T)$ is PI such that it is either a periodic orbit or  $X$ is an extension of an addition machine $X'$ such that $X'=X_{eq}$, 
    then $X\perp \mathcal{M}_{TM}$.
\item There is a weakly mixing extension of  an adding machine which is in $\mathcal{M}_{TM}^{\perp}$.  And there is a weakly mixing extension of  an adding machine which is not in $\mathcal{M}_{TM}^{\perp}$.
\end{enumerate}
\end{thm}
\begin{proof} (1) Suppose that $X\curlywedge \mathcal{M}_{TM}$. If $X_{eq}\neq X_{rat}$, then there is an irrational eigenvalue $\lambda$ of $(X,T)$. But $(\widehat{\Lambda_{\lambda}}, R_{\alpha})$ is a totally minimal system. Clearly, $X_{eq}\not\perp \widehat{\Lambda_{\lambda}}$. Then it follows from Lemma \ref{weaklydsijointmm}  that $X$ is not weakly disjoint from $\widehat{\Lambda_{\lambda}}$. This leads to a contradiction.

Now suppose that $X_{eq}=X_{rat}$. Then $X_{eq}$ is totally disconnected. But for any totally minimal system $(Y,S)$, $Y_{eq}$ is connected. Thus $X_{eq}\perp Y_{eq}$. By Lemma \ref{weaklydsijointmm}, one has $X\curlywedge Y$ and hence $X\curlywedge \mathcal{M}_{TM}$.




\medskip

(2) By (1) the collection of  $X\perp \mathcal{M}_{TM}$ is a subclass of the minimal systems with $X_{eq}=X_{rat}.$

If $X_{eq}$ is finite, then it is clear that $X\perp \mathcal{M}_{TM}$ if and only if $X$ is finite. To see this we note that if $X_{eq}$ is finite
and $X\not=X_{eq}$ then $X$ is a product of $X_{eq}$ with a non-trivial weakly mixing system $Z$ which implies that $X\not \perp Z$, and hence
$X\not\perp \mathcal{M}_{TM}.$

If $X_{eq}$ is not finite then it is an adding machine. Our assumption implies that $X_{eq}=X_{rat}$ and hence $X\curlywedge \mathcal{M}_{TM}$  by (1). So, it follows from Lemma \ref{weaklydsijointmm} (3)-(4) and the definition of $PI$ systems that $X\perp \mathcal{M}_{TM}$.



\medskip

(3) We need the construction by Glasner. Let $(Z,T)$ be an adding machine. By Glasner's result in \cite{G80}, there is a weakly mixing extension $(X,T)$ of $(Z,T)$ such that every minimal system disjoint from $(Z,T)$ is also disjoint from $(X,T)$. Since $Z\perp\mathcal{M}_{TM}$, we conclude that $Z\perp\mathcal{M}_{TM}$.

For the second assertion, we take the product of an adding machine and a weakly mixing system.



\end{proof}



In \cite[Theorem 6.3]{HSXY} it gives several equivalent conditions for a transitive system in $\mathcal{M}^\perp$.



\begin{thm} We have
\begin{enumerate}
\item Let $(X,T)$ be transitive and $X\perp \mathcal{M}$. Then $(X\times Y)\perp \mathcal{M}_{TM}$, where $(Y,S)$ is minimal and $Y\perp \mathcal{M}_{TM}$.
\item Let $(X,T)$ be scattering. Then $(X\times Y)\curlywedge \mathcal{M}_{TM}$, where $(Y,S)$ is minimal (or an $M$-system) and $Y\curlywedge \mathcal{M}_{TM}$.
\end{enumerate}
\end{thm}
\begin{proof} (1) Let $(Z,W)\in \mathcal{M}_{TM}$. Since $Y\perp Z$, $Y\times Z$ is minimal. Since $X\perp Y$, we know that for each transitive point $x\in X$ and any point $y\in Y$, $(x,y)$ is a transitive point of $X\times Y$. Moreover, and we know that $X\perp (Y\times Z)$.

Thus for any transitive point $(x,y)$ and any $z\in Z$,
$$\overline{\O}(((x,y),z), (T\times S)\times W)=\overline{\O}((x, (y,z)), T\times (S\times W))=X\times Y\times Z.$$
This show that $(X\times Y)\perp Z$.

\medskip


(2) Let $(Z,W)\in \mathcal{M}_{TM}$ and assume that $Y\curlywedge \mathcal{M}_{TM}$. It follows that the set of minimal points of $Y\times Z$ is dense in $Y\times Z$ and $Y\times Z$ is transitive, i.e. $Y\times Z$ is an $M$-system.
We claim that that $X\curlywedge (Y\times Z)$. To say the claim we note that for any open non-empty sets $U$ and $V$ of $Y\times Z$, there is a minimal set $M$ of $Y\times Z$ such that $M\cap U\not=\emptyset$ and $M\cap V\not=\emptyset$ by the transitivity of $Y\times Z$. So, by the fact that $X\curlywedge M$, we get that $X\curlywedge (Y\times Z)$. This ends the proof of the claim.

So, there is a transitive point $x\in X$ and a transitive point $(y,z)\in Y\times Z$ such that $(x,(y,z))$ is a transitive point of $X\times (Y\times Z)$. 
Then  we get
$$\overline{\O}(((x,y),z), (T\times S)\times W)=\overline{\O}((x, (y,z)), T\times (S\times W))=X\times Y\times Z.$$
This shows that $(X\times Y)\curlywedge  Z$.
\end{proof}

\section{Transitive systems disjoint from all totally minimal systems}
In this section we will give a characterization of transitive systems disjoint from all totally minimal systems or minimal weakly mixing systems. Let $\mathcal{P}$ to denote a dynamical property, e.g. weakly mixing or totally minimal. Let $\mathcal{M}_{\mathcal{P}}$ denote the collection of minimal systems with property $\mathcal{P}$. We will show the following theorem.

For an $M$-system, let $X_{\mathcal{P}}^{\perp}$ denote the closure of $\bigcup\{ M\subset X: \ M\text{ is a minimal set and }M\perp\mathcal{M}_\mathcal{P}\}$.
\begin{thm}\label{mainthm for dis}
Suppose that the property $\mathcal{P}$ satisfies that
\begin{itemize}
\item[(a)] every transitive system in $\mathcal{M}_{\mathcal{P}}^{\perp}$ is an $M$-system and
\item[(b)]  property $\mathcal{P}$ for minimal systems is preserved by almost one to one extension and by taking factors.
\end{itemize}
Then a transitive system $(X,T)$ is in $\mathcal{M}_{\mathcal{P}}^{\perp}$ if and only if
either $X$ has dense minimal sets $\{M_i\}_{i=1}^\infty$ such that $M_{i}\perp \mathcal{M}_{\mathcal{P}}$ for each $i\in \N$, or $X^\perp_{\mathcal{P}}$ is nowhere dense, and there are dense minimal sets $\{N_i\}_{i=1}^\infty$ such that $N_i\not\perp\mathcal{M}_{\mathcal{P}}$ and $X$ disjoint from all  quasifactors of $N_i$  having property $\mathcal{P}$ for each $i\in\mathbb{N}$.
\end{thm}

\begin{lem}
If $\mathcal{P}$ stands for weakly mixing or total minimality then  $\mathcal{P}$ satisfies conditions (a) and (b)
in Theorem \ref{mainthm for dis}.
\end{lem}
\begin{proof}
Let $\mathcal{P}$ stands for weakly mixing or total minimality. In \cite{HY05}, Huang and Ye show that for every thickly syndetic set $A\subset \N$, there is a subset $B\subset A$ such that ${\bf 1}_{B}$ is a minimal point in $\{0,1\}^{\mathbb{Z}_{+}}$. From this, they show that every transitive system in $\mathcal{M}^{\perp}$ is an $M$-system. However, $B$ constructed by Huang and Ye satisfies that $\overline{\O}({\bf 1}_{B},\sigma)$ is minimal weakly mixing, as pointed out in \cite[Remark 4.8]{DSY12}. In particular $\overline{\O}({\bf 1}_{B},\sigma)$ is also totally minimal. Thus  it follows from the proof of \cite[Theorem 2.6]{HY05} that a transitive system in $\mathcal{M}_{\mathcal{P}}^{\perp}$ is an $M$-system.
\end{proof}

In the sequel, we show Theorem \ref{mainthm for dis} for total minimality since it is similar for property $\mathcal{P}$ satisfying condition (a) and (b) in Theorem \ref{mainthm for dis}.

First we will consider the case when the transitive system is minimal, and then investigate the general case.

\subsection{Minimal case}

To consider the minimal case we need two known lemmas.

Let $(X,T)$ be minimal. Then any minimal subsystem of $(2^X,T)$ is called a quasi-factor of $X$.
\begin{lem}\cite[Chapter 11, Theorem 17]{Aus}\label{Glasner?} Let $(X,T)$ and $(Y, S)$ be minimal. If $X\not\perp Y$ then there is a HP extension $X^*$ of $X$ such that it has a factor which is a
non-trivial quasi factor of $Y$.
\end{lem}

\begin{lem}\cite[Chapter 11, Proposition 11]{Aus}\label{Glasner??} Let $(X,T)$ be minimal and $(N,T)$ be its non-trivial quasi factor. Then $X\not\perp N$.
\end{lem}

Note that non-trivial distal factor is equivalent to non-trivial equicontinuous factor. 
\begin{thm}\label{minimal-disjoint} A minimal system $(X,T)\in \mathcal{M}_{TM}^\perp$ 
if and only  if each non-trivial quasi-factor $(N,T)$ of $(X,T)$ is not totally  minimal (=$N_{eq}$ is not totally minimal).
\end{thm}
\begin{proof} ($\Leftarrow$) Assume that $(X,T)\not \perp (Y,S)$ for some $(Y,S)\in \mathcal{M}_{TM}$. Then by Lemma \ref{Glasner?} there is some HP extension $(Y^*,S)$ of $(Y,S)$ such that $(Y^*,S)$ has a factor $(N, S)$ which is a non-trivial quasi factor of $(X,T)$. Note that $(Y^*,S)\in \mathcal{M}_{TM}$ (as the extension is proximal, the maximal equicontinuous factor of $(Y^*,S)$ is the same as the one of $(Y,S)$). This implies that $(N,S)\in \mathcal{M}_{TM}$,
 a contradiction.

\medskip
($\Rightarrow$) Let $(X,T)\in \mathcal{M}_{TM}$. Assume the contrary that there is a non-trivial quasi factor $(N,T)$ which is totally minimal.
Then $(X,T)\perp (N,T)$, a contradiction since by Lemma \ref{Glasner??}, $(X,T)\not\perp (N,T)$. So, each non-trivial quasi factor of $X$ is
not totally minimal, i.e. $N_{eq}$ is not totally minimal.

\end{proof}

\begin{cor} Let $(X,T)$ be minimal and $X\not\in \mathcal{M}_{TM}$. Then there is a non-trivial quasi-factor of $X$ which is TM.
\end{cor}

\subsection{Transitive case and the proof of Theorem D}
\begin{lem} Let $(X,T)$ be transitive and $X\perp \mathcal{M}_{TM}$. Then $(X,T)$ is an $M$-system.
\end{lem}
\begin{proof}
By Lemma \ref{char of M-sys}, it suffices to show that $N(x,U)$ is piecewise syndetic for some $x\in Tran_{T}$ and any neighborhood $U$ of $x$. Now assume that $N(x,U)$ is not piecewise syndetic, then $A:=\N\setminus N(x,U)$ is thickly syndetic.  By \cite[Theorem 2.4]{HY05}, there is some nonempty subset $B\subset C$ such that ${\bf 1}_{B}$ is a minimal point in $\Sigma_2:=\{0,1\}^{\mathbb{Z}_{+}}$. Further, $\overline{\O}({\bf 1}_{B},\sigma)$ is minimal weakly mixing as pointed out in \cite[Remark 4.8]{DSY12}. But then $N({\bf 1}_{B}, V)=B\subset A$, where $V=\{(x_n)\in\Sigma_2:\ x_0=1\}$. This implies that $(U\times V)\cap  \O((x,{\bf 1}_{B}), T\times \sigma)=\emptyset$. Particularly, one has $X\perp \overline{\O}({\bf 1}_{B},\sigma)$ since  $X\perp \mathcal{M}_{TM}$ and $\overline{\O}({\bf 1}_{B},\sigma)$ is weakly mixing and hence totally minimal. This is a contradiction.
\end{proof}

Let $(X,T)$ be an $M$-system. Recall that we denote  by $X_{TM}^\perp$ the closure of
$$\bigcup\{M\subset X:\  M \text{ is minimal and}\ M\perp \mathcal{M}_{TM}\}. $$ 

\begin{thm}\label{main=thm}
Let $(X,T)$ be an $M$-system. Then the following statements are equivalent.
\begin{enumerate}
\item $(X,T)\perp \mathcal{M}_{TM}$.
\item Either $X$ has dense minimal sets $\{M_i\}_{i=1}^\infty$ such that $M_{i}\perp \mathcal{M}_{TM}$ for each $i\in \N$, or $X^\perp_{TM}$ is nowhere dense, and there are dense minimal sets $\{N_i\}_{i=1}^\infty$ such that $N_i\not\perp\mathcal{M}_{TM}$ and $X$ disjoint from all  totally minimal quasifactors of $N_i$ for each $i\in\mathbb{N}$.
\end{enumerate}
\end{thm}

\begin{proof}
 Let $MS_1:=\{M\perp \mathcal{M}_{TM}: M \text{ is a minimal set of } X\}$ and  $MS_2:=\{M\not\perp \mathcal{M}_{TM}: M \text{ is a minimal set of } X\}$.

(1) $\Longrightarrow (2)$.  We may assume that $\overline{\bigcup_{M\in MS_1}M}\neq X$. Then $\overline{\bigcup_{M\in MS_1}M}$ is nowhere dense. Otherwise, the interior of $\overline{\bigcup_{M\in MS_1}M}$ is not empty. Then there is some transitive point $x\in \overline{\bigcup_{M\in MS_1}M}$. Thus there is a sequence $(x_n)$ in $\bigcup_{M\in MS_1}M$ such that $x_n\rightarrow x$. But then $\overline{orb(x_n)}\rightarrow \overline{orb(x)}=X$ and hence $\overline{\bigcup_{M\in MS_1}M}=X$. This contradiction shows that $\overline{\bigcup_{M\in MS_1}M}$ is nowhere dense.
Further, we conclude that $\overline{\bigcup_{M\in MS_2}M}$ is dense.  Now that $X\perp \mathcal{M}_{TM}$, $X$ is disjoint from any totally minimal quasifactors of $N$ for each $N\in MS_2$.

(2) $\Longrightarrow (1)$. In the first case, suppose that there are dense minimal set $\{M_i\}$ such that $M_{i}\perp \mathcal{M}_{TM}$ for each $i$. Then it is clear that $X\perp \mathcal{M}_{TM}$. Indeed, let $Y$ be a totally minimal system and $J$ be a joining of $X$ and $Y$. Then $M_i\times Y\subset J$ for each $i$, since $M_i\perp Y$. Thus $X\times Y=J$ by the density of $\{M_i\}$. Thus $X\perp Y$ and hence $X\perp \mathcal{M}_{TM}$.

In the second case, suppose that there are dense minimal sets $\{N_i\}$ such that $N_i\not\perp\mathcal{M}_{TM}$ and $X$ disjoint from all  totally minimal quasifactors of $N_i$ for each $i\in\mathbb{N}$.  Then  it follows from Theorem \ref{key lemma} that $X\perp \mathcal{M}_{TM}$.
\end{proof}

\subsection{Quasifactors}

To make some notion simpler, we assume that systems consider in this section are homeomorphisms. This is not essential for study disjointness (see \cite{HY05}).
\begin{defn}Let $(X,T)$ be a minimal system and $\mathcal{X}$ is a subsystem of $2^{X}$. Let \{$\mathcal{X}_n\}_{i=1}^\infty$ be a collection of quasifactors of $X$. If for each $A\in \mathcal{X}$ there is some sequence $B_{n_i}\in \mathcal{X}_{n_i}$ such that $ \lim_{i\ra \infty} B_{n_i}\subset A$, then we say $\mathcal{X}\preceq \bigcup_{n=1}^{\infty}\mathcal{X}_{i}$.
\end{defn}

\begin{defn}
Let $(X,T)$ be a system. For $A\in 2^X$, define the {\em order of $A$} as follows:
$$\text{ord}(A)=\sup\ \{n\in \N: A\cap TA\cap \cdots \cap T^{n-1}A \not = \emptyset\}\in \N\cup\{\infty\}.$$
\end{defn}
Note that if $A$ is a minimal point in $2^X$ and is not trivial then $\text{ord}(A)<\infty$. Moreover, if $\mathcal{N}$ is a quasi-factor then $\text{ord}:\mathcal{N}\ra \N\cup\{\infty\},\ A\mapsto \text{ord}(A)$ is a constant function.

\medskip

The following lemma was proved in \cite[Lemma 5.5]{HSXY}.

\begin{lem}\label{small-disjoint}
Let $(X,T)$ be a minimal system. Suppose that  $\mathcal{X}$ is a subsystem of $(2^{X},T)$ and $\mathcal{Y}$ is a minimal subsystem of $(2^{X},T)$ (i.e. a quasifactor).  If $\mathcal{X}\perp \mathcal{Y}$, then for any $A\in \mathcal{X}$ and $B\in \mathcal{Y}$, we have
\[ A\cap  TB\cap \cdots\cap T^{{\rm ord}(\mathcal{Y})}B \neq\emptyset.\]
\end{lem}

\begin{lem}\label{lem-tran-1}
Let $(X,T)$ be a minimal system. Suppose that $\mathcal{X}$ and $\mathcal{Y}$ are transitive subsystems of $(2^{X}, T)$. If there exists a family of quasifactors $\{\mathcal{Y}_{i}\}_{i\in\mathbb{N}}$ such that $\mathcal{Y}\preceq \{\mathcal{Y}_i\}_{i\in\mathbb{N}}$, ${\rm ord}(\mathcal{Y}_i)\geq {\rm ord}(\mathcal{Y})$ and $\mathcal{X}\perp \mathcal{Y}_{i}$ for each $i\in\mathbb{N}$, then for any $A\in\mathcal{X}$ and $B\in\mathcal{Y}$ we have
\[ A\cap TB\cap \cdots\cap T^{{\rm ord}(\mathcal{Y})} B\neq\emptyset.\]
\end{lem}
\begin{proof}
Since $\mathcal{X}\perp \mathcal{Y}_i$ for each $i\in\N$, by Lemma \ref{small-disjoint} we have that for each $A\in \mathcal{X}$ and $C\in \mathcal{Y}_i$,
\[ A\cap TC\cap\cdots\cap T^{{\rm ord}(\mathcal{Y}_{i})}C\neq\emptyset.\]
In particular, since ${\rm ord}(\mathcal{Y}_i)\geq {\rm ord}(\mathcal{Y})$, it follows that
\begin{equation}\label{eq-1}
A\cap TC\cap\cdots\cap T^{{\rm ord}(\mathcal{Y})}C\neq\emptyset,
\end{equation}
for any $A\in \mathcal{X}$ and $C\in \bigcup_{i=1}^\infty \mathcal{Y}_i$.

Now since $\mathcal{Y}\preceq \{\mathcal{Y}_i\}_{i\in\mathbb{N}}$, for each $B\in \mathcal{Y}$ there is a sequence $C_{n_i}\in \mathcal{Y}_{n_i}$ with $\{n_i\}_{i\in \N}\subseteq \N$ such that $\lim_{i\to\infty}C_{n_i}\subset B$. Meanwhile, each $C_{n_i}$ satisfies (\ref{eq-1}). Taking limit as $i\to \infty$, we conclude that $A \cap TB \cap \cdots \cap T^{{\rm ord}(\mathcal{Y})}B \neq \emptyset$ for $A\in \mathcal{X}$ and $B\in \mathcal{Y}$. The proof is complete.
\end{proof}

\begin{thm}\label{lem-tran-2}
Let $(X,T)$ be a minimal system. Let $\{\mathcal{X}_\alpha\}_{\alpha\in \Gamma}$ be a collection of transitive subsystems of $(2^X, T)$. Assume that for each $\alpha\neq \beta\in \Gamma$, there exist a family of non-trivial quasifactors $\{\mathcal{Y}_i\}_{i\in \mathbb{N}}$ such that $\mathcal{X}_\beta\preceq \cup_{i\in\N}\mathcal{Y}_i$ in $2^X$, ${\rm ord}(\mathcal{Y}_i)\geq {\rm ord}(\mathcal{X}_\beta)$ and $\mathcal{X}_\alpha\perp \mathcal{Y}_i$ for each $i\in\mathbb{N}$.
Then $\Gamma$ is at most countable.
\end{thm}

\begin{proof}
 For $\alpha\in \Gamma$, we take a transitive point $A_{\alpha}$ of $(\mathcal{X}_{\alpha},T)$. Then by  the definition of ${\rm ord}(A_{\alpha})$, we have
$$\bigcap_{j=0}^{{\rm ord}(A_{\alpha})-1}T^jA_{\alpha}\neq \emptyset\text{ and } A_{\alpha}\cap TA_{\alpha}\cap \cdots \cap T^{{\rm ord}(A_{\alpha})}A_{\alpha}=\emptyset.$$ By Urysohn's Lemma, there is some $f_{\alpha}\in C(Y)$ such that
\[ f_{\alpha}\mid_{A_{\alpha}}\equiv 0\ \ \text{and}\ \ f_{\alpha}\mid_{TA_{\alpha}\cap \cdots \cap T^{{\rm ord}(A_{\alpha})}A_{\alpha}}\equiv 1.\]
Now, for any  $\alpha\neq\beta\in \Gamma$, by Lemma \ref{lem-tran-1}, we have
\[ \|f_{\alpha}-f_{\beta}\|_{\sup}\geq 1,\]
because $f_{\alpha}(x)=0$ and $f_{\beta}(x)=1$ for any $x\in  A_{\alpha}\cap TA_{\beta}\cap \cdots\cap T^{{\rm ord}(A_{\beta})}A_{\beta}$.
Thus, the separability of $C(X)$ implies the countability of $\Gamma$.
The proof is complete.
\end{proof}

\begin{lem}\label{yyyy}
Let $(X,T)$ be a non-trivial minimal system and $(\mathcal{W},T)$ be a non-trivial transitive subsystem of $(2^X,T)$. If $\{(\mathcal{X}_i,T)\}_{i\in\mathbb{Z}}$ is a sequence of quasifactors such that $\mathcal{W}\preceq\bigcup_{i=1}^\infty \mathcal{X}_i$ and ${\rm ord}(\mathcal{X}_i)\geq {\rm ord}(\mathcal{W})$ for each $i\in\mathbb{N}$, then
$\mathcal{W}\not\perp \mathcal{X}_{i_0}$ for some $i_0\in\mathbb{N}$.
\end{lem}

\begin{proof}
Otherwise, suppose that $\mathcal{W}\perp\mathcal{X}_i$ for each $i\in \mathbb{N}$.
Let $\Lambda$ be an uncountable index set. For each $\alpha\in \Lambda$, let $\mathcal{X}_{\alpha}={\mathcal W}$. This collection of non-trivial transitive subsystems $\{(\mathcal{X}_\alpha,T)\}_{\alpha\in \Lambda}$ of $(2^X,T)$ satisfies the conditions in Theorem \ref{lem-tran-2}. It follows that  $\Lambda$ is countable, a contradiction!
\end{proof}

Let $(X,T)$ be a transitive system and $M$ be a minimal subset of $X$. Let $(Y,T)$ be a minimal system. Suppose that $J$ is a nontrivial joining of $X$ and $Y$. Then let $\Omega$ and $\Omega_{M}$ be the sets of continuous points of the maps
\[ X\rightarrow 2^{Y},\ x\mapsto J[x] \text{ and } M\rightarrow 2^{Y},\ x\mapsto J[x],\]
respectively. Let
\[ W=\overline{\{(x, J[x])\in X\times 2^{Y}: x\in \Omega\}} \text{ and } W_{M}=\overline{\{(x, J[x])\in X\times 2^{Y}: x\in \Omega_{M}\}}.\]
 Then $\mathcal{Y}:=\pi_2(W)$ is a transitive subsystem of $2^{Y}$ and $\mathcal{Y}_{M}:=\pi_2(W_{M})$ is a quasifactor of $Y$.

 The following lemma is proved in Claim 1 of \cite[Proposition 5.10]{HSXY}.

 \begin{lem}\label{greater of order}
 ${\rm ord}(\mathcal{Y})\leq {\rm ord}(\mathcal{Y}_{M})$.
 \end{lem}

\begin{lem}\label{cover by order}
Let $(X,T)$ be a minimal system and $\mathcal{X}$ be a quasifactor of $X$. Then
\[ \bigcup_{A\in \mathcal{X}} \left(A\cap TA\cap \cdots\cap T^{{\rm ord}(\mathcal{X})-1}A\right)=X.\]
\end{lem}
\begin{proof}
We may assume that $\mathcal{X}$ is nontrivial and hence ${\rm ord}(\mathcal{X})=:k\in\mathbb{N}$.

Take $A\in\mathcal{X}$ and $x\in A\cap TA\cap \cdots\cap T^{k-1}A$. Now fix $y\in X$. Since $(X,T)$ is minimal, there is some sequence $(n_i)$ such that $T^{n_i}x\rightarrow y$. Without loss of generality, we may assume that $T^{n_i}A\rightarrow B$ in $2^{X}$. It remains to show that $y\in B\cap TB\cap \cdots\cap T^{k-1}B$. Note that
\[ \{x, T^{-1}x,\ldots, T^{-(k-1)}x\}\subset A.\]
Then for each $n\in \mathbb{Z}$,
\[ \{T^{n}x, T^{-1}T^{n}x,\ldots, T^{-(k-1)}T^{n}x\}\subset T^{n}A.\]
Thus
\begin{align*}
 \{y, T^{-1}y,\ldots,T^{-(k-1)}y\}&=\lim_{i\rightarrow\infty}\{T^{n_i}x, T^{-1}T^{n_i}x,\ldots, T^{-(k-1)}T^{n_i}x\}\\
 &\subset \lim_{i\rightarrow\infty} T^{n_i}A=B.
 \end{align*}
This implies that  $y\in B\cap TB\cap \cdots\cap T^{k-1}B$ and hence we complete the proof.
\end{proof}

\begin{lem}\label{equal of order}
Let $(X,T)$ and $(Y,T)$ be minimal systems. Let $J$ be a nontrivial joining of $X$ and $Y$. Then the canonically defined joining quasifactors $\mathcal{X}_{J}$ and $\mathcal{Y}_{J}$ of $X$ and $Y$ have the same order.
\end{lem}
\begin{proof}
Let $X_c$ and $Y_c$ be the set of continuous points of the maps
\[ X\rightarrow 2^{Y}, x\mapsto J[x] \text{ and } Y\rightarrow 2^{X}, y\mapsto J[y],\]
respectively. Let
\[ W=\overline{\{(x, J[x])\in X\times 2^{Y}: x\in X_c\}} \text{ and } V=\overline{\{(J[y],y)\in 2^X\times {Y}: y\in Y_{c}\}}.\]
Then $\mathcal{X}_{J}=\pi_1(V)$ and $\mathcal{Y}_{J}=\pi_2(W)$.
Let $k={\rm ord}(\mathcal{X}_J)$ and $\ell={\rm ord}(\mathcal{Y}_{J})$.

Now take $x\in X_c$ and then $ \ell={\rm ord}(J[x])$.
By Lemma \ref{cover by order}, there is some $A\in \mathcal{X}_{J}$ such that
\[ x\in A\cap TA\cap \cdots\cap T^{k-1}A.\]
By the construction of $V$, there is some $y\in Y$ such that $A\subset J[y]$. Thus
\[ x\in J[y]\cap TJ[y]\cap \cdots\cap T^{k-1} J[y].\]
This implies that
\[\left\{(x, y), (x, Ty),\ldots, (x,T^{k-1}y) \right\}\subset J.\]
Since $J$ is $T\times T$-invariant, we have
\[\left\{(T^{-(k-1)}x, y), (TT^{-(k-1)}x, y),\ldots, (T^{k-1}T^{-(k-1)}x, y) \right\}\subset J.\]
This implies that
\[ J[T^{-(k-1)}x]\cap TJ[T^{-(k-1)}x]\cap \cdots\cap T^{k-1}J[T^{-(k-1)}x]\neq\emptyset.\]
Thus ${\rm ord}(J[T^{-(k-1)}x] )\geq k$ and hence $\ell={\rm ord}(J[x])\geq k$.

\medskip

By the symmetry, we conclude that $k=\ell$.
\end{proof}

\subsection{The proof of the main theorem}

\begin{thm}\label{key lemma}
Let $(X,T)$ be an $M$-system. If there are dense minimal sets $\{N_i\}_{i=1}^\infty$ such that $N_i\not\perp\mathcal{M}_{TM}$ and $X$ is disjoint from all  totally minimal quasifactors of $N_i$ for each $i\in\mathbb{N}$, then $(X,T)\perp \mathcal{M}_{TM}$.
\end{thm}
\begin{proof}
To the contrary, assume that there is a totally minimal system $(Y,T)$ that is not disjoint from $(X,T)$.
Take a nontrivial joining $J\in Join(X, Y)$.

Let $\Omega$ be the set of continuous points of the map $X\rightarrow 2^{Y}, x\mapsto J[x]$. Let
\[ W:=\overline{\{(x, J[x])\in X\times 2^{Y}: x\in \Omega\}}.\]
Then $\pi_1: W\rightarrow X$ is almost one to one and $\mathcal{Y}:=\pi_2(W)$ is a nontrivial quasifactor of $Y$.

Since $\bigcup_{i=1}^\infty N_i$ is dense in $X$, by passing to some subsequence, we may assume that $J_i:=J\cap (N_i\times Y)$ is a nontrivial joining of $N_i$ and $Y$ for each $i\in\N$.

Let $X_{c,i}$ be the set of continuous points of the map $N_i\rightarrow 2^{Y}, x\mapsto J_i[x]=J[x]$ and let
\[ W_i:=\overline{\{(x,J[x])\in 2^{N_i}: x\in X_{c,i}\}}.\]
Let $\mathcal{Y}_{i}'=\pi_2(W_i)$.

Let $Y_{c,i}$ be the set of continuous points of the map $Y\rightarrow 2^{X}, y\mapsto J_i[y]=J[y]\cap N_i$ and let
\[ U_i:=\overline{\{(J_i[y],y)\in 2^{N_i}: y\in Y_{c,i}\}}.\]
Then $\mathcal{X}_i:=\pi_1(U_i)$ is a quasifactor of $N_i$. Since $\pi_2: U_i\rightarrow Y$ is almost one to one and $(Y,T)$ is totally minimal, we conclude that $\mathcal{X}_i$ is a non-trivial totally minimal quasifactor of $N_i$. Note that the map
\[\phi: \mathcal{X}_i\rightarrow 2^{Y}: A\mapsto \{ y\in Y: (A, y)\in U_i\}\]
is also upper semi-continuous. Thus the set $\mathcal{X}_{c,i}$ of continuous points of $\phi$ is residual in $\mathcal{X}_i$. Let
\[ V_i=\overline{\{(A, \phi_i(A))\in \mathcal{X}_i\times 2^{Y}: A\in \mathcal{X}_{c,i}\}}.\]
Then $\pi_1:V_i\rightarrow \mathcal{X}_i$ is almost one to one and $\mathcal{Y}_i:=\pi_2(V_i)$  is a quasifactor of $Y$.

\medskip

\noindent {\bf Claim 1}. $\mathcal{Y}\perp \mathcal{Y}_i$ for each $i\in\N$.
\begin{proof}[Proof of Claim 1]
Recall that $X\perp \mathcal{X}_i$ and hence $W\perp V_i$, since $W\rightarrow X$ and $V_i\rightarrow \mathcal{X}_i$ are both almost one to one. Further, $\mathcal{Y}\perp \mathcal{Y}_i$ since $\mathcal{Y}$ is a factor of $W$ and $\mathcal{Y}_i$ is a factor of $V_i$.
\end{proof}

\noindent {\bf Claim 2}. ${\rm ord}(\mathcal{Y}_i)\geq {\rm ord}(\mathcal{Y})$ for each $i\in\N$.
\begin{proof}[Proof of Claim 2]
By Lemma \ref{equal of order}, ${\rm ord}(\mathcal{X}_i)={\rm ord}(\mathcal{Y}_i')$. By Lemma \ref{greater of order}, ${\rm ord}(\mathcal{Y}_i')\geq {\rm ord}(\mathcal{Y})$. Hence  ${\rm ord}(\mathcal{X}_i)\geq {\rm ord}(\mathcal{Y})$. Thus it remains to show that ${\rm ord}(\mathcal{Y}_i)\geq {\rm ord}(\mathcal{X}_i)$.

Recall that $\pi_1: U_i\rightarrow \mathcal{X}_i$ is semi-open. Thus $\pi_1(Y_{c,i})$ is a residual set in $\mathcal{X}_i$. Thus there is some $y\in Y_{c,i}$ such that $J_i[y]\in \mathcal{X}_{c,i}$. Note that ${\rm ord}(\mathcal{X}_i)={\rm ord}(J_i[y])$. Thus
\[ J_i[y]\cap TJ_{i}[y]\cap \cdots\cap T^{{\rm ord}(\mathcal{X}_i)-1} J_{i}[y]\neq\emptyset.\]
This implies that
\begin{align*}
y&\in \phi_i(J_i[y])\cap \phi_i(TJ_i[y])\cap \cdots\cap \phi_i(T^{{\rm ord}(\mathcal{X}_i)-1}J_i[y])\\
&=\phi_i(J_i[y])\cap T\phi_i(J_i[y])\cap \cdots\cap T^{{\rm ord}(\mathcal{X}_i)-1}\phi_i(J_i[y]).
\end{align*}
In particular, $\phi_i(J_i[y])\cap T\phi_i(J_i[y])\cap \cdots\cap T^{{\rm ord}(\mathcal{X}_i)-1}\phi_i(J_i[y])\neq\emptyset$. Thus ${\rm ord}(\mathcal{Y}_i)\geq {\rm ord}(\mathcal{X}_i)$. This completes the proof of the claim.
\end{proof}

\noindent{\bf Claim 3}. $\mathcal{Y}\preceq \bigcup_{i\in\mathbb{N}} \mathcal{Y}_{i}$.
\begin{proof}[Proof of Claim 3]
For each $J_i[y]\in \mathcal{X}_{c,i}$ and $x\in J_i[y]$, we have $\phi_i(J_i[y])\subset J_i[x]$. Now for each $x\in \Omega$, take $J_i[y_{i_j}]\in \mathcal{X}_{c,i_j} $ and $x_{i_j}\in J_i[y_{i_j}]$ such that $x_{i_j}\rightarrow x$. Then $\lim \phi_{i_j} (J_{i_j}[y_{i_j}])\subset \lim J_{i_j}[x_{i_j}]\subset J[x]$. This shows that $\mathcal{Y}\preceq \bigcup_{i\in\mathbb{N}} \mathcal{Y}_{i}$.
\end{proof}

Now $\mathcal{Y}$ and $\{\mathcal{Y}_i\}$ satisfy the condition in Lemma \ref{yyyy}. But then the conclusion in Lemma \ref{yyyy} contradicts Claim 1. This contradiction leads to complete the proof.
\end{proof}

\begin{rem} There is a simple proof for the following statements: Assume $(X,T)$ is transitive. Then $X\perp \mathcal{M}$ if and only if
$(X,T)$ is an $M$-system and for any dense minimal sets $\{M_i\}$, $X\perp M_i$ (see \cite{HSXY}). We hope to find a simpler proof for Theorem \ref{main=thm}.
\end{rem}

\section{Summary of the results and counterparts in ergodic theory}

In this section we will summary the results in the previous sections and provide more questions for the further study.

\subsection{Summary of the results}

In this subsection we will summary of the results obtained in the previous sections.

\medskip

Summary of important sequences and families.
\medskip

\begin{center}
\begin{tabular}{|c|c|c|}
\hline
 Sequence & Family & Characterizations\\ \hline
 Birkhoff recurrence sequence & $\mathcal{F}_{Bir}$ &  Theorem \ref{equi of rec}   \\ \hline
 $R$-sequence & $\mathcal{F}_{R}$  &  Theorem \ref{rdop-char}  \\ \hline
 Dense orbit set & $\mathcal{F}_{DOS}$ & See \cite{HY05, GTWZ21, GL24,KRS22, XY22}\\ \hline
  $R$-sequence fo TM & $\mathcal{F}_{RTM}$& Theorem \ref{RTM-char} \\ \hline
\end{tabular}
\end{center}

\medskip
Let $(X,T)$ be a transitive system, $x\in Tran_{T}$ and $U$ be any nonempty open set of $X$.

\medskip

\begin{tabular}{|m{3cm}|m{5cm}|m{5cm}|}
\hline
 Disjoint classes & Characterization via recurrence & Other Characterizations\\ \hline
 $X\perp \mathcal{M}$ & $N(x, U)\in\mathcal{F}_{DOS}$ \cite{HY05} & See \cite{HSY20, HSXY}    \\ \hline
 $X\curlywedge \mathcal{M}$ & $N(x, U)\in\mathcal{F}_{R}$ (Theorem \ref{equi-s-point}) & Theorem \ref{scatterin-char-w} for $M$-systems  \\ \hline
 $X\perp \mathcal{M}_{TM}$ & $N(x, U)\cap A\neq\emptyset, \forall m$-set  $A$ for TM (Theorem \ref{rec char of dis TM}) & Theorem \ref{mainthm for dis} \\ \hline
 $X\curlywedge \mathcal{M}_{TM}$ & $N(x, U)\in\mathcal{F}_{RTM}$ (Theorem \ref{rec char of dis TM}) & Theorem \ref{weakdisjointtmminimal}  \\ \hline
\end{tabular}

\medskip
Note that while Katznelson and Richter's question remains open, stronger versions of the problem—which consider the question for an arbitrary minimal point $y$ in an arbitrary minimal system $(Y,S)$—have been resolved in the negative, as established in Proposition \ref{wprd-equ} and Corollary \ref{dis-equi-dis}.


\subsection{Counterparts in ergodic theory}
In this paper, we mainly focus on recurrence and orbits along  given sequences in topological dynamics. In a similar flavor, one may consider related questions in ergodic theory.

Recall that $A$ is a {\it Poincar\'e sequence} if for each measure preserving system (m.p.s. for short) $(X,\mathcal{B},\mu,T)$ and $B\in \mathcal{B}$ with $\mu(B)>0$ there is $n\in A$ such that $\mu(B\cap T^{-n}B)>0$. It is clear that a Poincar\'{e} sequence is a Birkhoff one. Similar to the Katznelson's question, Bergelson asked in  \cite{Berg85} that whether a Poincar\'e sequence for all ergodic m.p.s. with discrete spectrum is necessary a Poincar\'e sequence? However, this is answered negatively by Kriz in \cite{Kriz87}. For further study along this line, one may refer to Griesmer's work \cite{Gr21,Gr25}.

Similar to Richter's question, the ergodic counterpart is also far from being true. We start with some notions. $A\subset \N$ is called {\it a strongly Poincar\'e sequence} if $A\cap (n+(B-B))\not=\emptyset$ for any $B$ with positive upper Banach density and any $n\in \N$. It is clear that a strongly Poincar\'e sequence is an $R$-sequence. 

The following characterization is given in \cite[Lemma 6.1]{Gr19}.
\begin{thm} Let $A\subset \N$. Then the following statements are equivalent:
\begin{enumerate}
\item $A$ is a strongly Poincar\'e sequence.
\item $A$ is a shift invariant Poincar\'e sequence.
\item For each ergodic mps $(X,\mathcal{B},\mu,T)$, any $C, D\in \mathcal{B}$ with $\mu(C)\mu(D)>0$,
$$A\cap (\{k\in\N:\mu(C\cap T^{-k}D)>0\})\not=\emptyset.$$
\item For each ergodic m.p.s. $(X,\mathcal{B},\mu,T)$, any $B\in \mathcal{B}$ with $\mu(B)>0$ and any $n\in\N$,
$$A\cap (n+\{k\in\N:\mu(B\cap T^{-k}B)>0\})\not=\emptyset.$$
\item For each ergodic m.p.s. $(X,\mathcal{B},\mu,T)$, any $C\in \mathcal{B}$ with $\mu(C)>0$,
$$\mu(\cup_{a\in A}T^{-a}C)=1.$$
\end{enumerate}
\end{thm}

We want to know the connection between strongly Poincar\'e sequences and the ergodic averages.


\begin{prop}\label{ergodicarege} Assume that $A\subset \N$ satisfying that for any ergodic m.p.s.  $(X,\mathcal{B}, \mu,T)$ and any $f\in L^2$
$$\frac{1}{N}\sum_{i=1}^Nf(T^{a_i}x)\ra \int_X fd \mu \ \text{in}\ L^2.$$
Then $A$ is a strongly Poincar\'e sequence.
\end{prop}
\begin{proof} Assume that $\frac{1}{N}\sum_{i=1}^Nf(T^{a_i}x)\ra \int_X fd \mu.$ Then for any $g\in L^2$
$$\frac{1}{N}\sum_{i=1}^N\int_X gf(T^{a_i})d\mu\lra \int_X fd\mu \int_X gd\mu.$$

Now let $C,D\in \mathcal{B}$ with positive measures and $g=1_C$ and $f=1_D$. Then we have
$$\lim_{i\ra \infty} \frac{1}{N}\sum_{i=1}^N\mu(C\cap T^{-a_i}D)>0.$$
Thus, there is $i\in \N$ such that $\mu(C\cap T^{-a_i}D)>0.$ For a given $n\in\N$, set $D=T^nC$. Then we have that
$$A\cap (n+\{i\in\N: \mu(C\cap T^{-i}C)>0\})\not=\emptyset,$$
i.e. $A$ is a strongly Poincar\'e sequence.

\end{proof}

The converse of Proposition \ref{ergodicarege} is far from being true \cite{Gr19}.  However, it follows from \cite{ARS24} that  a strongly Poincar\'e sequence for any ergodic m.p.s. with discrete spectrum must be a strongly Poincar\'e sequence for all nilsystems. In view of Proposition \ref{simpleproof}, we may ask if a strongly Poincar\'e sequence for any ergodic m.p.s. with discrete spectrum must be a strongly Poincar\'e sequence for all measurable distal systems.

One may also refer to \cite{Ack22} for the generalization of Proposition \ref{ergodicarege} to countable abelian groups. 

\medskip


\begin{thebibliography}{SSS}
\bibitem{Ack22} E. M. Ackelsberg, {\it Rigidity, weak mixing, and recurrence in abelian groups}. Discrete Contin. Dyn. Syst. {\bf 42} (2022), no. 4, 1669-1705.

\bibitem{ARS24} E. M. Ackelsberg, F. K. Ritcher and O. Shalom, {\it On the maximal spectral type of nilsystems}. Proc. Amer. Math. Soc. Ser. B {\bf 11} (2024), 469-480.

\bibitem{Al25} R. Alweiss, {\it New obstacles to multiple recurrence}. arXiv: 2511.21680v1.


\bibitem{AG01} E. Akin and E. Glasner, {\it Residual properties and almost equicontinuity}. J. Anal. Math. {\bf 84} (2001), 243-286.

\bibitem{AG03} E. Akin and S. Kolyada, {\it  Li-Yorke sensitivity}. Nonlinearity, {\bf 16} (2003), no. 4, 1421-1433.

\bibitem{Aus} J. Auslander, {\it Minimal flows and their extensions}.  North-Holland mathematics studies; 153. (1988)


\bibitem{AF94} J. Auslander and H. Furstenberg, {\it Product recurrence and distal points}, Transactions of the American Mathematical Society,  {\bf 343} (1994), 221-232.

\bibitem{Berg85} B. Bergelson, {\it Open problem. Joint summer conference of AMS on Ramsey theory}. Arcata. 1985.

\bibitem{BR02} V. Bergelson and I. Rusza. {\it Squarefree numbers, IP sets, and ergodic theory}. In: Paul
Erd\H{s} and his Mathematics. Vol. 4. 5. Janos Bolyai Mathematical Society, 2002. Chap. 8, pp. 147-160.

\bibitem{BHM} F. Blanchard, B. Host and A. Maass, {\it  Topological complexity},  Ergodic Theory Dynam. Systems, {\bf 20} (2000), 641-662.


\bibitem{DSY12} P. Dong, S. Shao, and X. Ye. {\it Product recurrent properties, disjointness and weak disjointness}.
Israel J. Math., {\bf 188}(2012), 463–507.

\bibitem{DHM24} S. Donoso, F. Heren\'andez and A. Maass, {\it On recurrence for $\mathbb{Z}^d$ Weyl systems}, Monatshefte fur Mathematik {\bf 207}(2025):241–274.

\bibitem{F2}H. Furstenberg,  {\it Disjointness in ergodic theory, minimal sets, and a problem in Diophantine approximation,}
 Math. Systems Theory, {\bf 1} (1967), 1-49.

\bibitem{F1}
H. Furstenberg, \textit{Recurrence in ergodic theory and
combinatorial number theory}, M. B. Porter Lectures. Princeton
University Press, Princeton, N.J., 1981.

\bibitem{G80} S. Glasner, {\it Minimal skew products}. Trans. Amer. Math. Soc. {\bf 260} (1980), no. 2, 509-514.

\bibitem{GHSWY20}E. Glasner, W. Huang, S. Shao, and B. Weiss, X. Ye,
{\it Topological characteristic factors and nilsystems}, J. Eur. Math. Soc., {\bf 27} (2025), no. 1, 279-331.

\bibitem{GTWZ21} E. Glasner, T. Tsankov, B. Weiss and A. Zucker, {\it Bernoulli disjointness}, Duke Math. J., {\bf 170} (2021), no. 4,
 615-651.

\bibitem{GW93} E. Glasner and B. Weiss, {\it Sensitive dependence on initial conditions}. Nonlinearity {\bf 6} (1993), no. 6, 1067-1075.

\bibitem{GW17} E. Glasner and B. Weiss, {\it Three zutot}. Topol. Methods Nonlinear Anal. {\bf 49} (2017), no. 1, 351-358.

\bibitem{GW24} E. Glasner and  B. Weiss, \textit{On the class of systems which are disjoint from every ergodic system}, arXiv:2405.00463.

\bibitem{GKR22} D. Glasscock, A. Koutsogiannis and F. K. Richter, {\it On Katznelson’s Question for skew-product systems}, Bull. Amer. Math. Soc. (N.S.) {\bf 59}(2022), no. 4, 569-606.

\bibitem{GL24} D. Glasscock and Anh N. Le, {\it Dynamically syndetic sets and the combinatorics of syndetic, idempotent filters}, arXiv:2408.12785 v2.


\bibitem{Gr19} J. T.  Griesmer, {\it Recurrence, rigidity, and popular differences}. Ergodic Theory Dynam. Systems {\bf 39} (2019), no. 5, 1299-1316.

\bibitem{Gr21} J. T.  Griesmer,  {\it Separating Bohr denseness from measurable recurrence}. Discrete Anal. 2021, Paper No. 9, 20 pp.

\bibitem{Gr23} J. T.  Griesmer, {\it Special cases and equivalent forms of Katznelson's problem on recurrence}. Monatsh. Math. {\bf 200} (2023), no. 1, 63-79.

\bibitem{Gr25} J. T.  Griesmer, {\it Separating topological recurrence from measurable recurrence: exposition and extension of Kriz's example}. Australas. J. Combin. {\bf 91} (2025), 1-19.

\bibitem{GLR24} M. Gorska, M. Lemanczyk and T. de la Rue, {\it On orthogonality to uniquely ergodic systems},
 arXiv:2404.07907V3

\bibitem{HO08} K. Haddad and W. Ott, {\it Recurrence in pairs}, Ergodic Theory and Dynamical Systems, {\bf 28} (2008), 1135-1143.


\bibitem{HKM16} B. Host, B. Kra, and A. Maass, {\it Variations on topological recurrence}. Monatsh. Math., {\bf 179}(2016), 57-89.

\bibitem{HK18} B. Host and B. Kra, {\it Nilpotent structures in ergodic theory}. Mathematical Surveys and Monographs, Vol. 236 (2018).

 \bibitem{HSQY26} W. Huang, S. Shao, J. Qiu and X. Ye, {\it Nilpotent structures in dynamical systems},  Acta Math. Sin. (Engl. Ser.) to appear in 2026.

\bibitem{HSY04} W. Huang, S. Shao and X. Ye, {\it  Mixing and proximal cells along sequences}. Nonlinearity {\bf 17} (2004), no. 4, 1245-1260.

\bibitem{HSY16} W. Huang, S. Shao and X. Ye, {\it Nil Bohr-sets and almost automorphy of higher
order}, Mem. Amer. Math. Soc. {\bf 241} (2016), no. 1143, v+83 pp.

\bibitem{HSY20} W. Huang, S. Shao and X. Ye, {\it An answer to Furstenberg's problem on topological disjointness}. Ergodic Theory Dynam. Systems {\bf 40} (2020), 2467-2481.


\bibitem{HSXY} W. Huang, S. Shao, H. Xu and X. Ye, {\it On systems disjoint from all minimal systems}, arXiv:2504.17504.

\bibitem{HY02}  W. Huang and X. Ye, {\it An explicit scattering, non-weakly mixing example and weak disjointness}, Nonlinearity, {\bf 15}(2002), 849-862.

 \bibitem{HY04} W. Huang and X. Ye, {\it Topological complexity, return times and weak disjointness}, Ergodic Theory Dynam. Systems,  {\bf 24} (2004), no. 3, 825-846.

\bibitem{HY05} W. Huang and X. Ye, {\it Dynamical systems disjoint from any minimal system}. Trans. Amer. Math. Soc., {\bf 357}(2005),  669-694.

\bibitem{HY12} W. Huang and X. Ye, {\it Generic eigenvalues, generic factors and weak disjointness}. Dynamical systems and group actions, 119-142, Contemp. Math., 567, Amer. Math. Soc., Providence, RI, 2012.


\bibitem{Ka99} Y. Katznelson, {\it Chromatic numbers of Cayley graphs on $\Z$ and recurrence}. Com
binatorica, {\bf 21}(2001), 211-219. Paul Erdos and his mathematics (Budapest, 1999).

\bibitem{KRS22} M. Kennedy, S. Raum, and G. Salomon. {\it Amenability, proximality and higher-order syndeticity.} Forum Math. Sigma, 10:Paper No. e22, 28, 2022.

\bibitem{KTS01} S. Kolyada, S. Trofimchuk and L. Snoha, {\it Noinvertible minimal maps}. Fundamenta Mathematicae {\bf 168}(2): 141-163, 2001.

\bibitem{Kriz87} I. K\v{r}i\v{z}, {\it Large independent sets in shift-invariant graphs: solution of Bergelson’s problem}, Graphs Combin. {\bf 3} (1987), no. 2, 145-158

\bibitem{LYY15} J. Li, K. Yan and X. Ye, {\it Recurrence properties and disjointness on the induced spaces}, Discrete Contin. Dyn. Syst., {\bf 35} (2015), no. 3, 1059-1073.


\bibitem{MP13} G. L. Mullen and D. Panario, {\it Handbook of Finite Fields}. Chapman and Hall/CRC (2013).

\bibitem{Op19} P. Oprocha, {\it Double minimality, entropy and disjointness with all minimal systems}, Discrete Contin. Dyn. Syst., {\bf 39} (2019), no. 1, 263-275.


\bibitem{Peleg} R. Peleg, {\it Weak disjointness of transformation groups}. Proceedings of the American Mathematical Society, {\bf 33}(1972), No. 1, pp. 165- 170

\bibitem{QJ23}  J. Qiu,  {\it  Polynomial orbits for totally minimal systems}, Adv. Math., {\bf 432} (2023), 109260.

\bibitem{Veech70} W. A. Veech, {\it Point-distal flows}.  Amer. J. Math. {\bf 92} (1970), 205-242.

\bibitem{XY22} H. Xu and X. Ye,  {\it Disjointness from all minimal systems under group actions}, to appear in Israel J. of Math., arXiv:2212.07830.



\end{thebibliography}
\end{document}